\documentclass[11pt,a4paper,usenames,dvipsnames]{article}

\usepackage{url,amsmath,amssymb,latexsym,pstricks,mathrsfs,comment,amsthm,graphicx,tikz,tikz-cd,enumerate,accents,pgffor,cite,wrapfig,multicol,float,calc,geometry,tikz-3dplot}
\usepackage[colorlinks]{hyperref}
\usepackage[subnum]{cases}

\newcommand{\nc}{\newcommand}
\nc{\rnc}{\renewcommand}

\nc\U{{\raise1.6 ex\hbox{\rotatebox{180}{$U$}}}}

\usepackage{xcolor}

\usepackage{enumitem}

\usepackage[T1]{fontenc}
\usepackage{ wasysym }
\usepackage{stmaryrd}

\usepackage{tikz}
\usetikzlibrary{decorations.markings}
\usetikzlibrary{arrows,matrix}
\usepgflibrary{arrows}
  \usepgflibrary{fadings}
\usetikzlibrary{matrix,arrows,positioning,automata}
\pgfdeclarelayer{background layer}
\pgfsetlayers{background layer,main}

\tikzset{->-/.style={decoration={
  markings,
  mark=at position #1 with {\arrow{>}}},postaction={decorate}}}

\usetikzlibrary{intersections,through,hobby}

\usetikzlibrary{calc}

\allowdisplaybreaks

\nc\lamb{\overline\lam}
\nc\rhob{\overline\rho}

\rnc\L{\mathcal L}
\nc\R{\mathcal R}
\nc\bu{{\bf u}}
\nc\bv{{\bf v}}
\nc\bx{{\bf x}}
\nc\by{{\bf y}}
\nc\bz{{\bf z}}
\nc\bl{{\bf l}}
\nc\br{{\bf r}}
\nc\fnt{\lfloor\frac n2\rfloor}
\nc\s{\mathfrak s}
\rnc\t{\mathfrak t}
\newcommand{\hack}{\smash{\left\{\vphantom{\begin{aligned} \\[0.5ex] \\[0.5ex] \\[0.5ex] \end{aligned}}\right.}\!}
\nc\Rect{\mathscr R}

\nc\udotted[2]{\draw[dotted](#1+.4,2)--(#2-.4,2);}
\nc\ddotted[2]{\draw[dotted](#1+.4,0)--(#2-.4,0);}

\nc{\ubluebox}[2]{\bluebox{#1}{1.7}{#2}2\udotted{#1}{#2}}
\nc{\lbluebox}[2]{\bluebox{#1}0{#2}{.3}\ddotted{#1}{#2}}
\nc{\ublueboxes}[1]{{
\foreach \x/\y in {#1}
{ \ubluebox{\x}{\y}}}
}
\nc{\lblueboxes}[1]{{
\foreach \x/\y in {#1}
{ \lbluebox{\x}{\y}}}
}
\nc{\udottedsm}[2]{\draw [dotted] (#1+.4,2)--(#2-.4,2);}
\nc{\udottedsms}[1]{{
\foreach \x/\y in {#1}
{ \udottedsm{\x}{\y}}
}}
\nc{\ldottedsm}[2]{\draw [dotted] (#1+.4,0)--(#2-.4,0);}
\nc{\ldottedsms}[1]{{
\foreach \x/\y in {#1}
{ \ldottedsm{\x}{\y}}
}}
\nc{\bluebox}[4]{
\draw[color=blue!20, fill=blue!20] (#1,#2)--(#3,#2)--(#3,#4)--(#1,#4)--(#1,#2);
}
\nc{\uudotted}[2]{\draw [dotted] (#1+.4,4)--(#2-.4,4);}
\nc{\uubluebox}[2]{\bluebox{#1}{3.7}{#2}4\uudotted{#1}{#2}}
\nc{\uublueboxes}[1]{{
\foreach \x/\y in {#1}
{ \uubluebox{\x}{\y}}}
}
\nc{\uudottedsm}[2]{\draw [dotted] (#1+.4,4)--(#2-.4,4);}
\nc{\uudottedsms}[1]{{
\foreach \x/\y in {#1}
{ \uudottedsm{\x}{\y}}
}}

\nc\bn{{\bf n}}
\nc\Mod[1]{\ (\operatorname{mod}\ #1)}
\nc\TL{\mathcal T\!\mathcal L}
\newcommand{\darcxx}[4]{\draw[#4](#1,0)arc(180:90:#3) (#1+#3,#3)--(#2-#3,#3) (#2-#3,#3) arc(90:0:#3);}
\newcommand{\uarcxx}[4]{\draw[#4](#1,2)arc(180:270:#3) (#1+#3,2-#3)--(#2-#3,2-#3) (#2-#3,2-#3) arc(270:360:#3);}
\newcommand{\stlines}[1]{{\foreach \x/\y in {#1} { \stline{\x}{\y} }}}
\newcommand{\ustlines}[1]{{\foreach \x/\y in {#1} { \ustline{\x}{\y} }}}
\nc{\set}[2]{\{#1:#2\}}
\nc{\pres}[2]{\la#1:#2\ra}
\nc\sub\subseteq
\nc\mt\mapsto
\nc\al\alpha
\nc\be\beta
\nc\ga\gamma
\nc\de\delta
\nc\lam\lambda
\nc\Om\Omega
\nc\bit{\begin{itemize}}
\nc\eit{\end{itemize}}
\nc\ben{\begin{enumerate}[label=\textup{(\roman*)},leftmargin=7mm]}
\nc\bena{\begin{enumerate}[label=\textup{(\alph*)},leftmargin=7mm]}
\nc\een{\end{enumerate}}
\nc\bmc{\begin{multicols}}
\nc\emc{\end{multicols}}
\nc\bp{{\bf p}}
\rnc\iff{\ \Leftrightarrow\ }
\rnc\implies{\ \Rightarrow\ }
\nc\pf{\begin{proof}}
\nc\epf{\end{proof}}
\nc\epfres{\hfill\qed}
\let\oldproofname=\proofname
\renewcommand{\proofname}{\rm\bf{\oldproofname}}
\nc\AND{\qquad\text{and}\qquad}
\nc{\COMMA}{,\qquad}
\nc\ol\overline
\nc\wh\widehat
\rnc\emptyset\varnothing
\nc\la\langle
\nc\ra\rangle

\newcommand{\uv}[1]{\fill (#1,2)circle(.17);}
\newcommand{\uuv}[1]{\fill (#1,4)circle(.17);}
\newcommand{\lv}[1]{\fill (#1,0)circle(.17);}
\newcommand{\uvs}[1]{{\foreach \x in {#1} { \uv{\x}}}}
\newcommand{\uuvs}[1]{{\foreach \x in {#1} { \uuv{\x}}}}
\newcommand{\lvs}[1]{{\foreach \x in {#1} { \lv{\x}}}}
\newcommand{\darcx}[3]{\draw(#1,0)arc(180:90:#3) (#1+#3,#3)--(#2-#3,#3) (#2-#3,#3) arc(90:0:#3);}
\newcommand{\darc}[2]{\darcx{#1}{#2}{.4}}
\newcommand{\uarcx}[3]{\draw(#1,2)arc(180:270:#3) (#1+#3,2-#3)--(#2-#3,2-#3) (#2-#3,2-#3) arc(270:360:#3);}
\newcommand{\uuarcx}[3]{\draw(#1,4)arc(180:270:#3) (#1+#3,4-#3)--(#2-#3,4-#3) (#2-#3,4-#3) arc(270:360:#3);}
\newcommand{\uarc}[2]{\uarcx{#1}{#2}{.4}}
\newcommand{\stline}[2]{\draw(#1,2)--(#2,0);}
\newcommand{\ustline}[2]{\draw(#1,4)--(#2,0);}
\newcommand{\uustline}[2]{\draw(#1,4)--(#2,0);\uuv{#1}}

\newcommand{\dom}{\operatorname{dom}} 
\newcommand{\codom}{\operatorname{codom}}
\newcommand{\rank}{\operatorname{rank}}
\newcommand\id{\operatorname{id}}

\newcommand{\sm}{\setminus}

\numberwithin{equation}{section}

\newtheorem{thm}[equation]{Theorem}
\newtheorem{lemma}[equation]{Lemma}
\newtheorem{cor}[equation]{Corollary}
\newtheorem{prop}[equation]{Proposition}

\theoremstyle{definition}

\newtheorem{rem}[equation]{Remark}

\begin{document}

\title{\vspace{-.5cm}Presentations for Temperley-Lieb algebras}
\author{James East\\
{\it\small Centre for Research in Mathematics and Data Science,}\\
{\it\small Western Sydney University, Locked Bag 1797, Penrith NSW 2751, Australia.}\\
{\tt\small J.East\,@\,WesternSydney.edu.au}}
\date{}

\maketitle
\begin{abstract}
We give a new and conceptually straightforward proof of the well-known presentation for the Temperley-Lieb algebra, via an alternative new presentation.  Our method involves twisted semigroup algebras, and we make use of two apparently new submonoids of the Temperley-Lieb monoid.



\end{abstract}

\thanks{Dedicated to the memory of Prof.~Vaughan F.~R.~Jones.}


\section{Introduction}

Temperley-Lieb algebras were introduced by their namesakes in \cite{TL1971} to study lattice models in statistical mechanics.  These algebras have since appeared naturally and extensively in many mathematical contexts, and have found particularly strong applications in knot theory; see especially the works of Kauffman and Jones, such as \cite{Jones1994_a,Jones1987,Jones1983_2,Jones1994, Kauffman1990,Kauffman1987,Kauffman1997}.  More information may be found in surveys such as \cite{Abramsky2008,RSA2014,KL1994}; the introductions to \cite{Westbury1995,BDP2002} also contain valuable discussions.

One of the most important tools for working with Temperley-Lieb algebras is a well-known presentation by generators and relations (stated in Theorem \ref{t:TL3} below), which provides an algebraic axiomatisation of the topological definition in terms of homotopy-classes of planar tangles.  
As far as the author is aware, the only complete proof of this presentation is the one given by Borisavljevi{\'c}, Do{\v{s}}en and Petri{\'c} in \cite{BDP2002}, where they also discuss the lack of a proof up to that point, even (surprisingly) of the fact that the so-called \emph{hooks} or \emph{diapsides} form a generating set; see the tangles denoted $\ol e_i$ in Figure~\ref{f:lamb}.  The proof given in \cite{BDP2002} involves working with a (countably infinite) monoid of tangles in place of the algebra itself, and is rather ingeneous, albeit quite involved, as the authors themselves note at the end of their introduction: \emph{``Our proof exhibits some difficulties, which we think cannot be evaded''}.

The main purpose of the current article is to give an alternative proof of the presentation, which we believe exhibits no such difficulties.  Our method also involves replacing the algebra with a monoid, but this time with a finite one, the so-called \emph{Temperley-Lieb monoid}, $\TL_n$.  This monoid is sometimes called the \emph{Jones monoid} in the literature, and denoted $\mathcal J_n$; see for example \cite{DEEFHHLM2019,Auinger2014,EG2017,EMRT2018}, and especially \cite{LF2006} for a discussion of naming conventions.  
The main innovation in our proof is in the use of two apparently new submonoids~$\L_n$ and~$\R_n$ of $\TL_n$; roughly speaking, these each capture half of the complexity of $\TL_n$ itself, and we have a natural factorisation $\TL_n=\L_n\R_n$ (Proposition \ref{p:TL_LR}).  We first give presentations for $\L_n$ and $\R_n$ in Section \ref{s:LR} (Theorem \ref{t:LR}), and in Section \ref{s:TL1} show how to stitch these together to obtain a new presentation for $\TL_n$ (Theorem \ref{t:TL1}).  In Section \ref{s:TL2}, we show how to rewrite the new presentation to obtain the original one (Theorem \ref{t:TL2}), and in Section \ref{ss:statements} explain how to convert any presentation for the monoid into a presentation for the algebra, using general results on twisted semigroup algebras \cite{JEgrpm}.

Some aspects of the proof presented here bear similarities to that given in \cite{BDP2002}, the most obvious being the normal forms used; compare our Lemma \ref{l:wxy} with \cite[Lemmas 1 and 4--6]{BDP2002}.  Although these forms utilise different generating sets, they both capture natural combinatorial data associated to the strings of a tangle, as encoded in certain tuples defined in Section \ref{ss:tuples}; however, our normal forms have length at most $n$, while those of \cite{BDP2002} grow quadratically in $n$.  
As an application of our results, we are also able to establish the (well-known) equivalence of the topological framework of tangles and the combinatorial approach via set partitions; compare our Corollary \ref{c:TL_LR} with \cite[Remark 3]{BDP2002}.
Beyond these parallels, however, the author believes that the current proof is rather simpler, and hopes it may therefore be of some benefit.

\section{Preliminaries and statement of the main results}\label{s:prelim}

\subsection{Temperley-Lieb monoids and algebras}\label{ss:TL}

For the duration of the article, we fix a positive integer $n$, and write $\bn=\{1,\ldots,n\}$.  To avoid trivialities we assume that $n\geq3$.

For $i\in\bn$, we define the points in the plane by $P_i=(i,1)$ and~$P_i'=(i,0)$.  
A \emph{string} is a smooth non-self-intersecting embedding $\s$ of the unit interval $[0,1]$ into the rectangle $\Rect_n=[1,n]\times[0,1]$ such that 
\bit
\item $\s(0),\s(1)\in\{P_1,\ldots,P_n,P_1',\ldots,P_n'\}$, and 
\item other than these endpoints, the image of $\s$ is contained in the interior of $\Rect_n$.  
\eit
A \emph{planar tangle} is a collection $\al=\{\s_1,\ldots,\s_n\}$ of $n$ pairwise non-intersecting strings.  Such a planar tangle is typically identified with the union of the images of its strings, which is a collection of~$n$ curves contained in $\Rect_n$.  Note that $\{\s_1(0),\ldots,\s_n(0),\s_1(1),\ldots,\s_n(1)\}=\{P_1,\ldots,P_n,P_1',\ldots,P_n'\}$.  Some examples with $n=9$ are given in Figure \ref{f:albe}.

Two planar tangles $\al=\{\s_1,\ldots,\s_n\}$ and $\be=\{\t_1,\ldots,\t_n\}$ are \emph{equivalent}, written as $\al\equiv\be$, if there is a homotopy through planar tangles from $\al$ to $\be$, or more formally if there is a family of continuous maps $F_i:[0,1]\times[0,1]\to\Rect_n$, $i\in\bn$, such that, relabelling the strings if necessary, and writing $\s_i^u(x)=F_i(u,x)$ for $i\in\bn$ and $u,x\in[0,1]$,
\bit
\item $\ga_u = (\s_1^u,\ldots,\s_n^u)$ is a planar tangle for all $u\in[0,1]$,
\item $\ga_0=\al$ and $\ga_1=\be$.
\eit
We denote by $\TL_n$ the set of all $\equiv$-classes of planar tangles.

The product $\al\be$ of two planar tangles $\al$ and $\be$ is formed by shifting $\al$ one unit in the positive $y$-direction, attaching it to $\be$, rescaling to the whole object until it lies in the rectangle $\Rect_n$, and finally removing any loops contained entirely in the interior of $\Rect_n$.  (For later use, let $m(\al,\be)$ be the number of such loops.)  Figure \ref{f:albe} gives an example calculation with $n=9$ (and where $m(\al,\be)=1$).  It is clear that $\al_1\be_1\equiv\al_2\be_2$ whenever $\al_1\equiv\al_2$ and $\be_1\equiv\be_2$.  Thus, we have an induced product on $\equiv$-classes, which is easily seen to be associative.  It follows that $\TL_n$ is a semigroup under this product, indeed a monoid whose identity is the $\equiv$-class of the planar tangle $\id_n$ with $n$ vertical strings.  This monoid is called the \emph{Temperley-Lieb monoid (of degree $n$)}.  In what follows, we typically identify a planar tangle with its $\equiv$-class.

\begin{figure}[ht]
\begin{center}
\begin{tikzpicture}[scale=.43]
\begin{scope}[shift={(0,0)}]	
\uvs{1,...,9}
\lvs{1,...,9}
\stline13
\stline86
\stline99
\uarcx27{.6}
\uarcx34{.3}
\uarcx56{.3}
\darc12
\darc45
\darc78
\draw(0.6,1)node[left]{$\al=$};
\draw[->](10.5,-1)--(12.5,-1);
\end{scope}
\begin{scope}[shift={(0,-4)}]	
\uvs{1,...,9}
\lvs{1,...,9}
\uarc12
\uarc34
\uarc56
\uarc89
\darc12
\darcx45{.3}
\darcx36{.6}
\darc89
\stline77
\draw(0.6,1)node[left]{$\be=$};
\end{scope}
\begin{scope}[shift={(13,-1)}]	
\uvs{1,...,9}
\lvs{1,...,9}
\stline13
\stline86
\stline99
\uarcx27{.6}
\uarcx34{.3}
\uarcx56{.3}
\darc12
\darc45
\darc78
\draw[->](10.5,0)--(12.5,0);
\end{scope}
\begin{scope}[shift={(13,-3)}]	
\uvs{1,...,9}
\lvs{1,...,9}
\uarc12
\uarc34
\uarc56
\uarc89
\darc12
\darcx45{.3}
\darcx36{.6}
\darc89
\stline77
\end{scope}
\begin{scope}[shift={(26,-1)},yscale=0.5]	
\stline13
\stline86
\stline99
\uarcx27{.6}
\uarcx34{.3}
\uarcx56{.3}
\darc45
\darc78
\end{scope}
\begin{scope}[shift={(26,-2)},yscale=0.5]	
\uarc34
\uarc56
\uarc89
\darc12
\darcx45{.3}
\darcx36{.6}
\darc89
\stline77
\end{scope}
\begin{scope}[shift={(26,-2)}]	
\uvs{1,...,9}
\lvs{1,...,9}
\draw(9.4,1)node[right]{$=\al\be$};
\end{scope}
\end{tikzpicture}
\caption{Calculating a product $\al\be$, where $\al,\be\in\TL_9$.}
\label{f:albe}
\end{center}
\end{figure}

A \emph{transversal} of a planar tangle $\al$ is a string connecting points $P_i$ and $P_j'$ for some $i,j\in\bn$.  Every other string of~$\al$ is an \emph{upper} or \emph{lower non-transversal}, with obvious meanings.  We write $\rank(\al)$ for the number of transversals, noting that $\rank(\al)\equiv n\Mod2$.  We also define the domain and codomain,
\begin{align*}
\dom(\al) &= \set{i\in\bn}{\text{$P_i$ is the endpoint of a transversal of $\al$}},\\
\codom(\al) &= \set{i\in\bn}{\text{$P_i'$ is the endpoint of a transversal of $\al$}}.
\end{align*}
Note that $\rank(\al)=|{\dom(\al)}|=|{\codom(\al)}|$.  With $\al\in\TL_9$ as in Figure \ref{f:albe}, we have $\dom(\al)=\{1,8,9\}$, $\codom(\al)=\{3,6,9\}$ and $\rank(\al)=3$.
Certain obvious identities will come in handy, including
\[
\dom(\al\be)\sub\dom(\al) \COMMA \codom(\al\be)\sub\codom(\be) \COMMA \rank(\al\be)\leq\min(\rank(\al),\rank(\be)).
\]
Another is that $\codom(\al)\sub\dom(\be)\implies\rank(\al\be)=\rank(\al)$.  These will typically be used without explicit reference.  

There is a natural involution $\TL_n\to\TL_n:\al\mt\al^\dagger$, defined by reflection in the line $y=\frac12$.  See Figure \ref{f:*} for an example when $n=9$.  Since $(\al^\dagger)^\dagger=\al=\al\al^\dagger\al$ and $(\al\be)^\dagger=\be^\dagger\al^\dagger$ for all $\al,\be$, it follows that $\TL_n$ is a \emph{regular $*$-semigroup} in the sense of \cite{NS1978}.  We also have identities such as
\[
 \dom(\al^\dagger)=\codom(\al) \COMMA
 \codom(\al^\dagger)=\dom(\al) \COMMA
 \rank(\al^\dagger)=\rank(\al).
\]
The symmetry/duality afforded by the involution will allow us to simplify many proofs.

\begin{figure}[ht]
\begin{center}
\begin{tikzpicture}[scale=.43]
\begin{scope}[shift={(0,0)}]	
\uvs{1,...,9}
\lvs{1,...,9}
\stline13
\stline86
\stline99
\uarcx27{.6}
\uarcx34{.3}
\uarcx56{.3}
\darc12
\darc45
\darc78
\draw(0.6,1)node[left]{$\al=$};
\draw[->](10.5,1)--(12.5,1);
\end{scope}
\begin{scope}[shift={(13,0)}]	
\uvs{1,...,9}
\lvs{1,...,9}
\stline31
\stline68
\stline99
\darcx27{.6}
\darcx34{.3}
\darcx56{.3}
\uarc12
\uarc45
\uarc78
\draw(9.4,1)node[right]{$=\al^\dagger$};
\end{scope}
\end{tikzpicture}
\caption{The involution $\al\mt\al^\dagger$.}
\label{f:*}
\end{center}
\end{figure}

Now let $\Bbbk$ be a field, and $\de$ an arbitrary element of $\Bbbk$.  The \emph{Temperley-Lieb algebra} $\TL_n(\Bbbk,\de)$ is the vector space over $\Bbbk$ with basis $\TL_n$, and with product $\star$ defined on basis elements $\al,\be\in\TL_n$ (and then extended by $\Bbbk$-linearity) by
\[
\al\star\be = \de^{m(\al,\be)} (\al\be),
\]
where $m(\al,\be)$ was defined above.  So the product in $\TL_n(\Bbbk,\de)$ of two basis elements is always a scalar multiple of a basis element, meaning that $\TL_n(\Bbbk,\de)$ is a so-called \emph{twisted semigroup algebra} \cite{Wilcox2007}.

\subsection{Presentations}\label{ss:pres}

For a set $X$, we write $X^*$ for the free monoid over $X$, which consists of all words over $X$ under concatenation; the empty word will be denoted by $1$.  For a set $\Om\sub X^*\times X^*$ of pairs of words, we write $\Om^\sharp$ for the congruence on $X^*$ generated by $\Om$.  We say a monoid $M$ has presentation $\pres X\Om$ if $M\cong X^*/\Om^\sharp$: i.e., if there is a monoid surmorphism $X^*\to M$ with kernel $\Om^\sharp$.  If $\phi$ is such a surmorphism, we say $M$ has presentation $\pres X\Om$ via $\phi$.  Elements of $X$ and $\Om$ are called generators and relations, respectively, and a relation $(u,v)\in\Om$ will often be denoted as an equation: $u=v$.

If $\Bbbk$ is a field, then the free (unital) $\Bbbk$-algebra over $X$ is $\Bbbk[X^*]$, the semigroup algebra of $X^*$.  For ${\Om\sub\Bbbk[X^*]\times\Bbbk[X^*]}$, we say an algebra $A$ has presentation $\pres X\Om$ if there is an algebra surmorphism $\phi:\Bbbk[X^*]\to A$ with kernel (in the vector space sense) spanned by $\set{u\phi-v\phi}{(u,v)\in\Om}$.

\subsection{Statement of the main results}\label{ss:statements}

The main results of this paper are presentations for the Temperley-Lieb monoid (Theorems \ref{t:TL1} and \ref{t:TL2}) and algebra (Theorem \ref{t:TL3}).  Theorem \ref{t:TL1} is new, and will be used to prove the well-known Theorems~\ref{t:TL2} and \ref{t:TL3}.

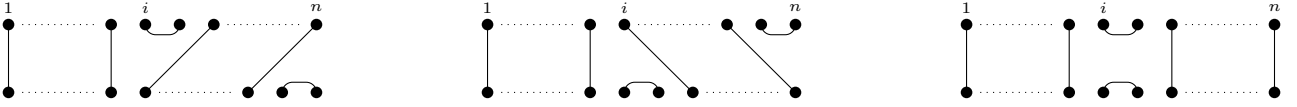
\begin{figure}[ht]
\begin{center}
\begin{tikzpicture}[scale=.45]
\begin{scope}[shift={(0,0)}]	
\uvs{1,4,5,6,7,10}
\lvs{1,4,5,8,9,10}
\udotted14
\ddotted14
\udotted7{10}
\ddotted58
\stline11
\stline44
\stline75
\stline{10}8
\uarcx56{.3}
\darcx9{10}{.3}
\node()at(1,2.5){\tiny$1$};
\node()at(5,2.5){\tiny$i$};
\node()at(10,2.5){\tiny$n$};
\end{scope}
\begin{scope}[shift={(14,0)}]	
\uvs{1,4,5,8,9,10}
\lvs{1,4,5,6,7,10}
\udotted14
\ddotted14
\ddotted7{10}
\udotted58
\stline11
\stline44
\stline57
\stline8{10}
\darcx56{.3}
\uarcx9{10}{.3}
\node()at(1,2.5){\tiny$1$};
\node()at(5,2.5){\tiny$i$};
\node()at(10,2.5){\tiny$n$};
\end{scope}
\begin{scope}[shift={(28,0)}]	
\uvs{1,4,5,6,7,10}
\lvs{1,4,5,6,7,10}
\udotted14
\ddotted14
\udotted7{10}
\ddotted7{10}
\stline11
\stline44
\stline77
\stline{10}{10}
\uarcx56{.3}
\darcx56{.3}
\node()at(1,2.5){\tiny$1$};
\node()at(5,2.5){\tiny$i$};
\node()at(10,2.5){\tiny$n$};
\end{scope}
\end{tikzpicture}
\caption{The planar tangles $\lamb_i$ (left), $\rhob_i$ (middle) and $\ol e_i$ (right), for $1\leq i\leq n-1$.}
\label{f:lamb}
\end{center}
\end{figure}

To state these results, we begin by defining three alphabets:
\[
L = \{\lam_1,\ldots,\lam_{n-1}\} \COMMA R = \{\rho_1,\ldots,\rho_{n-1}\} \COMMA E = \{e_1,\ldots,e_{n-1}\}.
\]
We think of the elements of these sets as abstract letters, but to each $1\leq i\leq n-1$, we associate three basic planar tangles, $\lamb_i$, $\rhob_i$ and $\ol e_i$, from $\TL_n$, as shown in Figure \ref{f:lamb}.

It is easy to see that
\[
\lamb_i^\dagger=\rhob_i \COMMA \rhob_i^\dagger=\lamb_i \COMMA \ol e_i^\dagger=\ol e_i \COMMA
\lamb_i\rhob_i = \ol e_i \COMMA \rhob_i\lamb_i = \ol e_{n-1} \qquad\text{for all $i$.}
\]
We define two morphisms
\[
\phi:(L\cup R)^*\to\TL_n:\lam_i\mt\lamb_i,\ \rho_i\mt\rhob_i \AND \psi:E^*\to\TL_n:e_i\mt\ol e_i,
\]
and we extend the over-line notation to words, writing $\ol u=u\phi$ for $u\in(L\cup R)^*$, and $\ol v=v\psi$ for $v\in E^*$.

Now consider the set $\Om$ consisting of the following relations over $L\cup R$, where in each relation the subscripts range over all meaningful values, subject to any stated constraints:
\begin{alignat}{2}
\tag*{(L1)} \label{L1} \lam_i\lam_{n-1} &= \lam_i \\ 
\tag*{(L2)} \label{L2} \lam_i\lam_j &= \lam_{j+2}\lam_i &&\qquad\qquad\text{if $i\leq j\leq n-3$}\\
\tag*{(L3)} \label{L3} \lam_{n-2i+1}^i\lam_{n-2i} &= \lam_{n-2i+1}^i \\[4truemm]
\tag*{(R1)} \label{R1} \rho_{n-1}\rho_i &= \rho_i \\ 
\tag*{(R2)} \label{R2} \rho_j\rho_i &= \rho_i\rho_{j+2} &&\qquad\qquad\text{if $i\leq j\leq n-3$}\\
\tag*{(R3)} \label{R3} \rho_{n-2i}\rho_{n-2i+1}^i &= \rho_{n-2i+1}^i \\[4truemm]
\tag*{(RL1)}\label{RL1}     &\mathrel{\phantom{=}}\hphantom{\hack} \lam_{n-1}\lam_j\rho_{i-2} 	&&\qquad\qquad\text{if $j\leq i-2$}\\ 
\tag*{(RL2)}\label{RL2} \rho_i\lam_j &=\hack \lam_{n-1}=\rho_{n-1} &&\qquad\qquad\text{if $i-1\leq j\leq i+1$} \\
\tag*{(RL3)}\label{RL3}     &\mathrel{\phantom{=}}\hphantom{\hack} \lam_{n-1}\lam_{j-2}\rho_i &&\qquad\qquad\text{if $j\geq i+2$.}
\end{alignat}			
Our main new result is the following, expressed in terms of the above notation:

\begin{thm}\label{t:TL1}
The Temperley-Lieb monoid $\TL_n$ has presentation $\pres{L\cup R}\Om$ via $\phi$.
\end{thm}

Now consider the set $\Xi$ consisting of the following relations over $E$:
\begin{align}
\tag*{(E1)} \label{E1} e_i^2 &= e_i &&\hspace{-2cm}\text{for all $i$}\\
\tag*{(E2)} \label{E2} e_ie_j &= e_je_i &&\hspace{-2cm}\text{if $|i-j|>1$}\\
\tag*{(E3)} \label{E3} e_ie_je_i &= e_i &&\hspace{-2cm}\text{if $|i-j|=1$.}
\end{align}

\begin{thm}\label{t:TL2}
The Temperley-Lieb monoid $\TL_n$ has presentation $\pres{E}\Xi$ via $\psi$.
\end{thm}

From any presentation for the monoid $\TL_n$, it is easy to deduce a presentation for the algebra $\TL_n(\Bbbk,\de)$, for any field $\Bbbk$ and any non-zero scalar $\de\in\Bbbk\sm\{0\}$, using results of \cite[Section~6]{JEgrpm}.  We define
\[
\Psi:\Bbbk[E^*]\to\TL_n(\Bbbk,\de)
\]
to be the $\Bbbk$-linear extension of $\psi$.  By \cite[Theorem 44]{JEgrpm}, the monoid presentation $\pres E\Xi$ for $\TL_n$ may be transformed into an algebra presentation $\pres E{\Xi'}$ for $\TL_n(\Bbbk,\de)$ by replacing each relation $u=v$ from $\Xi$ by $\de^{m(v)}u=\de^{m(u)}v$; here for a word $w=e_{i_1}\cdots e_{i_k}\in E^*$, $m(w)$ is defined to be the number of loops created when forming the product $\ol e_{i_1}\cdots\ol e_{i_k}$ in $\TL_n$.  The only relation from $\Xi$ that needs to be modified in this way is \ref{E1}, so let $\Xi'$ be the set of relations obtained from $\Xi$ by replacing~\ref{E1} with
\begin{equation}
\tag*{(E1)$'$} \label{E1'} e_i^2 = \de e_i \qquad\qquad\qquad\text{for all $i$.}
\end{equation}

\begin{thm}\label{t:TL3}
For $\de\in\Bbbk\sm\{0\}$, the Temperley-Lieb algebra $\TL_n(\Bbbk,\de)$ has presentation $\pres{E}{\Xi'}$ via $\Psi$.  \epfres
\end{thm}

\subsection{Partitions and tuples}\label{ss:tuples}

Certain combinatorial data associated to planar tangles will be of use throughout.  First, we write $\bn'=\{1',\ldots,n'\}$.  To each string $\s$ of $\al\in\TL_n$, we associate a two-element subset $\bp(\s)$ of $\bn\cup\bn'$: 
\bit
\item If $\s$ joins $P_i$ to $P_j$ for some $i,j\in\bn$, then $\bp(\s)=\{i,j\}$.
\item If $\s$ joins $P_i$ to $P_j'$ for some $i,j\in\bn$, then $\bp(\s)=\{i,j'\}$.
\item If $\s$ joins $P_i'$ to $P_j'$ for some $i,j\in\bn$, then $\bp(\s)=\{i',j'\}$.
\eit
We then define $\bp(\al)=\set{\bp(\s)}{\s\in\al}$; this is clearly a (set) partition of $\bn\cup\bn'$, whose blocks all have size $2$; cf.~\cite{Brauer1937}.  It turns out that the partition $\bp(\al)$ provides enough information to uniquely specify $\al$ itself, as follows from Corollary \ref{c:TL_LR} below.  

In fact, this partition provides more information than is strictly necessary.
Let the left-most endpoints of the upper and lower non-transversals of $\al\in\TL_n$ be $P_{x_1},\ldots,P_{x_k}$ and $P_{y_1}',\ldots,P_{y_k}'$, where $x_1>\cdots>x_k$ and $y_1>\cdots>y_k$.  We then define $\bl(\al)=(x_1,\ldots,x_k)$ and ${\br(\al)=(y_1,\ldots,y_k)}$.  By convention, if $k=0$ (i.e., if $\al$ has no non-transversals), we write ${\bl(\al)=\br(\al)=\emptyset}$.  Note that $\bl(\al^\dagger)=\br(\al)$ and $\br(\al^\dagger)=\bl(\al)$.

For example, with $\al,\be\in\TL_9$ as in Figure \ref{f:albe}, we have
\[
\bl(\al)=(5,3,2) \COMMA \br(\al)=(7,4,1) \COMMA \bl(\be)=(8,5,3,1) \COMMA \br(\be)=(8,4,3,1).
\]
As a foreshadowing of things to come, the reader may check that $\al=\lamb_5\lamb_3\lamb_2\rhob_1\rhob_4\rhob_7$ (cf.~Lemma~\ref{l:wxy}).

It will be important to know which tuples can occur as $\bl(\al)$ or $\br(\al)$ for some $\al\in\TL_n$.  With this in mind, let $T_n$ be the set of integer tuples $\bx=(x_1,\ldots,x_k)$ such that
\[
k\geq0 \COMMA  x_1>\cdots>x_k\geq1 \COMMA x_i\leq n-2i+1\text{ for all $i$.}
\]
The third item says that $x_1\leq n-1$, $x_2\leq n-3$, $x_3\leq n-5$, and so on.  For such a tuple $\bx$, we write $|\bx|=k$.  Note that $0\leq |\bx| \leq \fnt$.

\begin{lemma}\label{l:Tn}
If $\al\in\TL_n$, then $\bl(\al)$ and $\br(\al)$ both belong to $T_n$.
\end{lemma}

\pf
We just prove the statement for $\bl(\al)$, as that for $\br(\al)$ is dual.  Write $\bl(\al)=(x_1,\ldots,x_k)$.  If $k=0$ then $\bl(\al)=\emptyset\in T_n$, so for the rest of the proof we assume that $k\geq1$.  By definition, we have $x_1>\cdots>x_k\geq1$, and for each $1\leq i\leq k$, $\al$ has a string joining $P_{x_i}$ to $P_{y_i}$ for some $x_i<y_i\leq n$.  It quickly follows that $x_i=\min\{x_1,\ldots,x_i,y_1,\ldots,y_i\}$.  For $A\sub\bn$ with $|A|=2i$, the greatest possible value $\min(A)$ could take is $n-2i+1$, so the result follows.
\epf

Conversely, it follows from Corollary \ref{c:*} below that for any $\bx,\by\in T_n$ with $|\bx|=|\by|$, there exists $\al\in\TL_n$ with $\bl(\al)=\bx$ and $\br(\al)=\by$ (and from Corollary \ref{c:TL_LR} that this $\al$ is unique).

\section{The monoids $\L_n$ and $\R_n$}\label{s:LR}

In this section and the next we prove Theorem \ref{t:TL1}, which gives a presentation for the Temperley-Lieb monoid~$\TL_n$.  Our strategy involves two submonoids $\L_n$ and $\R_n$ of $\TL_n$; these will be defined shortly, and the main goal of the current section is to give presentations for them (Theorem~\ref{t:LR}).  These presentations will then be used in the next section, together with a factorisation ${\TL_n=\L_n\R_n}$ (Proposition~\ref{p:TL_LR}), to complete the proof of Theorem~\ref{t:TL1}.  The approach just described is based on that of \cite{JEinsn}, which employed similar methods to treat monoids of (planar/order-preserving) partial bijections.

We begin by defining the submonoids $\L_n$ and $\R_n$.  We say a planar tangle $\al\in\TL_n$ is:
\bit
\item \emph{right-simple} if $\br(\al)=(n-2k+1,n-2k+3,\ldots,n-3,n-1)$ for some $k\geq0$,
\item \emph{left-simple} if $\bl(\al)=(n-2k+1,n-2k+3,\ldots,n-3,n-1)$ for some $k\geq0$.
\eit
Equivalently, $\al\in\TL_n$ is right-simple if $\codom(\al)=\{1,\ldots,r\}$ for some $r\equiv n\Mod2$, and the partition $\bp(\al)$ contains the blocks ${\{r+1,r+2\}',\{r+3,r+4\}',\ldots,\{n-1,n\}'}$.  Similar comments apply to left-simple planar tangles.  Note that $\al$ is right-simple if and only if $\al^\dagger$ is left-simple, and vice versa.  Also note that $\lamb_i$ is right-simple, and $\rhob_i$ left-simple, for any $i$.

We define $\L_n$ (respectively, $\R_n$) to be the set of all right-simple (respectively, left-simple) planar tangles.  
Figure \ref{f:gade} shows elements $\ga\in\L_9$ and $\de\in\R_9$; note that $\ga\de$ is equal to the planar tangle $\al\in\TL_9$ pictured in Figure \ref{f:albe}.

\begin{figure}[ht]
\begin{center}
\begin{tikzpicture}[scale=.43]
\begin{scope}[shift={(0,0)}]	
\uuvs{1,...,9}
\lvs{1,...,9}
\ustline11
\ustline82
\ustline93
\uuarcx27{.6}
\uuarcx34{.3}
\uuarcx56{.3}
\darc45
\darc67
\darc89
\draw(0.6,2)node[left]{$\ga=$};
\end{scope}
\begin{scope}[shift={(16,0)}]	
\uuvs{1,...,9}
\lvs{1,...,9}
\ustline13
\ustline26
\ustline39
\uuarcx45{.4}
\uuarcx67{.4}
\uuarcx89{.4}
\darc12
\darc45
\darc78
\draw(0.6,2)node[left]{$\de=$};
\end{scope}
\end{tikzpicture}
\caption{Right- and left-simple planar tangles $\ga\in\L_9$ and $\de\in\R_9$.}
\label{f:gade}
\end{center}
\end{figure}

\begin{lemma}\label{l:LR}
\ben
\item \label{l:LR1} If $\al\in\L_n$, and if $\bl(\al)=(x_1,\ldots,x_k)$, then $\al=\lamb_{x_1}\cdots\lamb_{x_k}$.
\item \label{l:LR2} If $\al\in\R_n$, and if $\br(\al)=(x_1,\ldots,x_k)$, then $\al=\rhob_{x_k}\cdots\rhob_{x_1}$.
\een
\end{lemma}

\pf
We just prove \ref{l:LR1}, as \ref{l:LR2} is dual.  The proof is by induction on $k$.  If $k=0$, then $\al$ has no non-transversals, so we may write $\al=(\s_1,\ldots,\s_n)$, where for some permutation $\pi$ of $\bn$, the string $\s_i$ joins $P_i$ to $P_{i\pi}'$ for each $i\in\bn$.  If $\pi$ was not the identity permutation, say with $i<j$ and $i\pi>j\pi$, then the strings $\s_i$ and $\s_j$ would cross, a contradiction.  So $\pi$ is the identity permutation, and by planarity $\al=\id_n$.  The result holds in this case, since the stated product is empty.

Now suppose $k\geq1$.  Let $\s$ be the string of $\al$ with $P_{x_1}$ as one of its endpoints, and let the other endpoint of $\s$ be $P_y$.  We first claim that $y=x_1+1$.  To see this, note first that $y\geq x_1+1$ by definition of~$\bl(\al)$.  Aiming for a contradiction, suppose $y\geq x_1+2$, and let $\t$ be the string of $\al$ with $P_{x_1+1}$ as one of its endpoints.  Since $x_1<x_1+1<y$, and since $\s$ joins $P_{x_1}$ to $P_y$, it follows by planarity that $\t$ is a non-transversal, and is contained entirely in the region bounded by $\s$ and the upper border of the rectangle $\Rect_n$.  In particular, $\t$ joins $P_{x_1+1}$ to $P_z$ for some $x_1+1<z<y$.  But then $x_1+1$ must be one of the entries in the tuple $\bl(\al)$, and this contradicts the maximality of $x_1$.

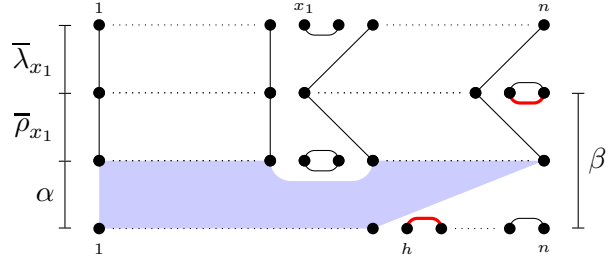
\begin{figure}[ht]
\begin{center}
\begin{tikzpicture}[scale=.45]
\begin{scope}[shift={(0,0)}]	
\darcxx{10}{11}{.3}{very thick,red}
\fill[blue!20](1,0)--(9,0)--(14,2)--(9,2)arc(360:270:.6)--(6+.6,2-.6)arc(270:180:.6)--(1,2)--(1,0);
\uvs{1,6,7,8,9,14}
\lvs{1,9,10,11,13,14}
\udotted16
\udotted9{14}
\ddotted19
\ddotted{11}{13}
\uarcx78{.3}
\darcx{13}{14}{.3}
\draw(0,1)node[left]{$\al$};
\draw[|-|](0,2)--(0,0);
\draw(15,2)node[right]{$\be$};
\draw[|-|](15,4)--(15,0);
\node()at(1,-.6){\tiny$1$};
\node()at(10,-.6){\tiny$h$};
\node()at(14,-.6){\tiny$n$};
\end{scope}
\begin{scope}[shift={(0,2)}]	
\uarcxx{13}{14}{.3}{very thick,red}
\uvs{1,6,7,12,13,14}
\lvs{1,6,7,8,9,14}
\udotted16
\udotted7{12}
\ddotted9{14}
\ddotted16
\darcx78{.3}
\stline11
\stline66
\stline79
\stline{12}{14}
\draw(0,1)node[left]{$\rhob_{x_1}$};
\draw[|-](0,2)--(0,0);
\end{scope}
\begin{scope}[shift={(0,4)}]	
\lvs{1,6,7,12,13,14}
\uvs{1,6,7,8,9,14}
\ddotted16
\ddotted7{12}
\udotted9{14}
\udotted16
\uarcx78{.3}
\darcx{13}{14}{.3}
\stline11
\stline66
\stline97
\stline{14}{12}
\draw(0,1)node[left]{$\lamb_{x_1}$};
\draw[|-](0,2)--(0,0);
\node()at(1,2.5){\tiny$1$};
\node()at(7,2.5){\tiny$x_1$};
\node()at(14,2.5){\tiny$n$};
\end{scope}
\end{tikzpicture}
\caption{Partitions constructed during the proof of Lemma \ref{l:LR}.  For simplicity, we write $h=n-2k+1$.}
\label{f:LR}
\end{center}
\end{figure}

So we have shown that $P_{x_1}$ is joined to $P_{x_1+1}$ by a string of $\al$.  Because of this string, it follows that $\al=\ol e_{x_1}\al = (\lamb_{x_1}\rhob_{x_1})\al$; cf.~Figure \ref{f:LR}.  Let $\be=\rhob_{x_1}\al$, and note that $\be$ contains a string joining $P_{n-1}$ to~$P_n$, and a string joining $P_{n-2k+1}'$ to $P_{n-2k+2}'$; these strings are coloured red in Figure \ref{f:LR}.  In fact, by maximality of $x_1$, we have
\[
\bl(\be)=(n-1,x_2,\ldots,x_k) \AND \br(\be)=\br(\al)=(n-2k+1,n-2k+3,\ldots,n-3,n-1).
\]
(Again, see Figure \ref{f:LR}.)  Let $\ga$ be the planar tangle obtained from $\be$ by replacing the two strings just mentioned by two transversals: one joining $P_{n-1}$ to $P_{n-2k+1}'$, and the other joining $P_n$ to $P_{n-2k+2}'$.  
Since $\lamb_{x_1}$ has a string joining $P_{n-1}'$ to $P_n'$, it follows that $\lamb_{x_1}\be=\lamb_{x_1}\ga$.  Moreover, we have 
\[
\bl(\ga)=(x_2,\ldots,x_k) \AND \br(\ga)=(n-2k+3,\ldots,n-3,n-1).
\]
The latter gives $\ga\in\L_n$.  By induction, since $\bl(\ga)$ has length $k-1$, we have $\ga=\lamb_{x_2}\cdots\lamb_{x_k}$.  It follows that $\al = (\lamb_{x_1}\rhob_{x_1})\al = \lamb_{x_1}\be = \lamb_{x_1}\ga = \lamb_{x_1}\cdot\lamb_{x_2}\cdots\lamb_{x_k}$, and the proof is complete.
\epf

\begin{cor}\label{c:LR}
\ben
\item \label{c:LR1} If $\al,\be\in\L_n$, then $\al=\be\iff\bl(\al)=\bl(\be)$.
\item \label{c:LR2} If $\al,\be\in\R_n$, then $\al=\be\iff\br(\al)=\br(\be)$.
\een
\end{cor}

\pf
We just prove \ref{c:LR1}, as \ref{c:LR2} is dual.  The forwards implication is clear, while the converse follows from Lemma \ref{l:LR}\ref{l:LR1}: if $\bl(\al)=\bl(\be)=(x_1,\ldots,x_k)$, then $\al=\lamb_{x_1}\cdots\lamb_{x_k}=\be$.
\epf

In what follows, we will write $\ol L=\{\lamb_1,\ldots,\lamb_{n-1}\}$ and $\ol R=\{\rhob_1,\ldots,\rhob_{n-1}\}$.

\begin{prop}\label{p:L}
The sets $\L_n$ and $\R_n$ are submonoids of $\TL_n$, and we have $\L_n=\la\ol L\ra$ and $\R_n=\la\ol R\ra$.
\end{prop}

\pf
As usual, it suffices to prove the statement for $\L_n$.  In light of Lemma \ref{l:LR}, it remains to show that $\la\ol L\ra\sub\L_n$.  Since $\ol L\sub\L_n$, it is enough to show that $\L_n$ is closed under right multiplication by elements of $\ol L$.  With this in mind, let $\al\in\L_n$, and let $1\leq i\leq n-1$ be arbitrary; we must show that $\al\lamb_i\in\L_n$.  Also let $r=\rank(\al)$.  If $r=n$ then $\al=\id_n$ (see the first paragraph of the proof of Lemma~\ref{l:LR}), and so $\al\lamb_i=\lamb_i\in\L_n$.  For the rest of the proof, we assume that $r<n$, and we note that
\[
\codom(\al) = \{1,\ldots,r\} \AND \br(\al) = (r+1,r+3,\ldots,n-3,n-1).
\]
Using this, we show in Figure \ref{f:L} that $\al\lamb_i\in\L_n$ in the cases:
\bena\bmc2
\item \label{La} $i<r$,
\item \label{Lb} $i=r$,
\item \label{Lc} $i>r$ and $i\not\equiv n\Mod2$,
\item \label{Ld} $i>r$ and $i\equiv n\Mod2$,
\emc\een
respectively.  In fact, in all but case \ref{La}, we have $\al\lamb_i=\al$.
\epf

\begin{figure}[ht]
\begin{center}
\begin{tikzpicture}[scale=.4]
\begin{scope}[shift={(-22,0)}]	
\fill[blue!20](1,0)--(10,0)--(16,4)--(1,4)--(1,0);
\uustline{2.5}1
\uustline54
\uustline{6.5}5
\uustline{8.5}6
\uustline{10}7
\uustline{13}{10}
\uuvs{1,16}
\lvs{1,4,5,6,7,10,11,12,15,16}
\uudotted1{16}
\darcx{11}{12}{.3}
\darcx{15}{16}{.3}
\draw(0,2)node[left]{$\al$};
\draw[|-](0,4)--(0,0);
\node()at(1,4.5){\tiny$1$};
\node()at(16,4.5){\tiny$n$};
\node()at(0,5){\ref{La}};
\end{scope}
\begin{scope}[shift={(-22,-2)}]	
\uvs{1,4,5,6,7,10,11,12,15,16}
\lvs{1,4,5,8,9,10,13,14,15,16,16}
\udotted14
\udotted7{10}
\udotted{12}{15}
\ddotted14
\ddotted58
\ddotted{10}{13}
\uarcx56{.3}
\stlines{1/1,4/4,7/5,10/8,11/9,12/10,15/13,16/14}
\darcx{15}{16}{.3}
\draw(0,1)node[left]{$\lamb_i$};
\draw[|-|](0,2)--(0,0);
\end{scope}
\begin{scope}[shift={(-22,-10)}]	
\fill[blue!20](1,0)--(8,0)--(16,4)--(1,4)--(1,0);
\uustline41
\uustline{9.5}7
\uustline{12}8
\uuvs{1,16}
\lvs{1,7,8,9,10,11,12,15,16}
\uudotted1{16}
\darcx9{10}{.3}
\darcx{11}{12}{.3}
\darcx{15}{16}{.3}
\draw(0,2)node[left]{$\al$};
\draw[|-](0,4)--(0,0);
\node()at(1,4.5){\tiny$1$};
\node()at(16,4.5){\tiny$n$};
\node()at(0,5){\ref{Lb}};
\end{scope}
\begin{scope}[shift={(-22,-12)}]	
\uvs{1,7,8,9,10,11,12,15,16}
\lvs{1,7,8,9,10,13,14,15,16}
\udotted17
\udotted{12}{15}
\ddotted17
\ddotted{10}{13}
\uarcx89{.3}
\stlines{1/1,7/7,10/8,11/9,12/10,15/13,16/14}
\darcx{15}{16}{.3}
\draw(0,1)node[left]{$\lamb_i$};
\draw[|-|](0,2)--(0,0);
\end{scope}
\begin{scope}[shift={(0,0)}]	
\fill[blue!20](1,0)--(4,0)--(18,4)--(1,4)--(1,0);
\uustline51
\uustline{10}4
\uuvs{1,18}
\lvs{1,4,5,6,9,10,11,12,13,14,17,18}
\uudotted1{18}
\darcx56{.3}
\darcx9{10}{.3}
\darcx{11}{12}{.3}
\darcx{13}{14}{.3}
\darcx{17}{18}{.3}
\draw(0,2)node[left]{$\al$};
\draw[|-](0,4)--(0,0);
\node()at(1,4.5){\tiny$1$};
\node()at(18,4.5){\tiny$n$};
\node()at(0,5){\ref{Lc}};
\end{scope}
\begin{scope}[shift={(0,-2)}]	
\uvs{1,4,5,6,9,10,11,12,13,14,17,18}
\lvs{1,4,5,6,9,10,11,12,15,16,17,18}
\udotted14
\udotted69
\udotted{14}{17}
\ddotted14
\ddotted69
\ddotted{12}{15}
\stlines{1/1,4/4,5/5,6/6,9/9,10/10,13/11,14/12,17/15,18/16}
\uarcx{11}{12}{.3}
\darcx{17}{18}{.3}
\draw(0,1)node[left]{$\lamb_i$};
\draw[|-|](0,2)--(0,0);
\end{scope}
\begin{scope}[shift={(0,-10)}]	
\fill[blue!20](1,0)--(4,0)--(18,4)--(1,4)--(1,0);
\uustline51
\uustline{10}4
\uuvs{1,18}
\lvs{1,4,5,6,9,10,11,12,13,14,17,18}
\uudotted1{18}
\darcx56{.3}
\darcx9{10}{.3}
\darcx{11}{12}{.3}
\darcx{13}{14}{.3}
\darcx{17}{18}{.3}
\draw(0,2)node[left]{$\al$};
\draw[|-](0,4)--(0,0);
\node()at(1,4.5){\tiny$1$};
\node()at(18,4.5){\tiny$n$};
\node()at(0,5){\ref{Ld}};
\end{scope}
\begin{scope}[shift={(0,-12)}]	
\uvs{1,4,5,6,9,10,11,12,13,14,17,18}
\lvs{1,4,5,6,9,10,11,12,15,16,17,18}
\udotted14
\udotted69
\udotted{14}{17}
\ddotted14
\ddotted69
\ddotted{12}{15}
\stlines{1/1,4/4,5/5,6/6,9/9,12/10,13/11,14/12,17/15,18/16}
\uarcx{10}{11}{.3}
\darcx{17}{18}{.3}
\draw(0,1)node[left]{$\lamb_i$};
\draw[|-|](0,2)--(0,0);
\end{scope}
\end{tikzpicture}
\caption{Verification that $\al\lamb_i\in\L_n$, as in the proof of Proposition \ref{p:L}.}
\label{f:L}
\end{center}
\end{figure}

Now that we know $\L_n$ and $\R_n$ are monoids, we wish to establish presentations.  Keeping Proposition~\ref{p:L} in mind, we have surmorphisms
\[
\phi_L:L^*\to\L_n:\lam_i\mt\lamb_i \AND \phi_R:R^*\to\R_n:\rho_i\mt\rhob_i,
\]
which are restrictions of the morphism $\phi:(L\cup R)^*\to\TL_n$ defined in Section \ref{ss:statements}.  Let $\Om_L$ and $\Om_R$ be the sets of relations \ref{L1}--\ref{L3} and \ref{R1}--\ref{R3}, respectively.  Our main goal in this section is to prove the following:

\begin{thm}\label{t:LR}
The monoids $\L_n$ and $\R_n$ have presentations $\pres L{\Om_L}$ and $\pres R{\Om_R}$, via $\phi_L$ and $\phi_R$, respectively.
\end{thm}

We will just prove Theorem \ref{t:LR} for $\L_n$, as the result for $\R_n=\L_n^\dagger$ is dual.  Since $\phi_L$ is a surmorphism, it remains to show that $\ker(\phi_L)=\Om_L^\sharp$.  For the rest of this section, we write ${\sim_L}=\Om_L^\sharp$.  

\begin{lemma}\label{l:relL}
We have ${\sim_L}\sub\ker(\phi_L)$.
\end{lemma}

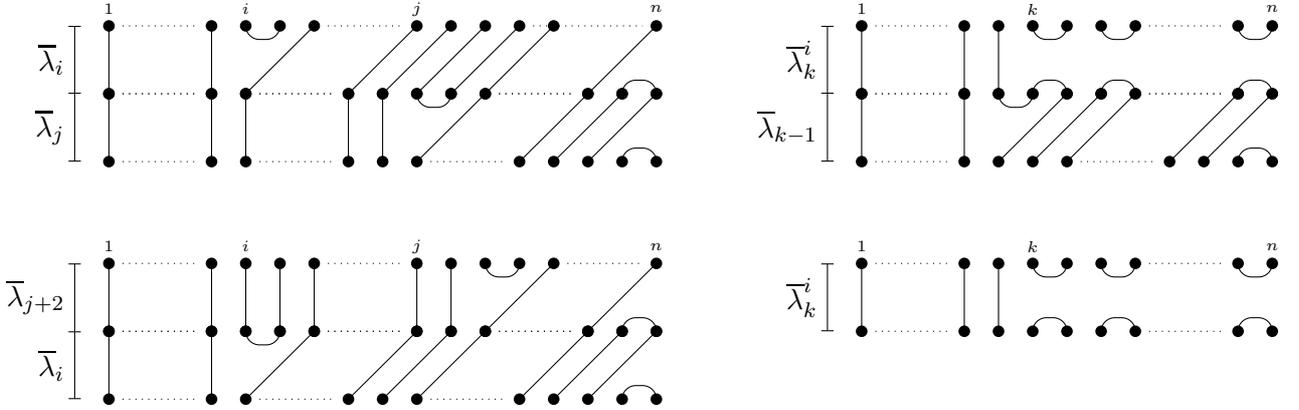
\begin{figure}[ht]
\begin{center}
\begin{tikzpicture}[scale=.45]
\begin{scope}[shift={(0,0)}]	
\uvs{0,3,4,5,6,9,10,11,12,13,16}
\lvs{0,3,4,7,8,9,10,11,14,15,16}
\udotted03
\udotted69
\udotted{12}{16}
\ddotted03
\ddotted47
\ddotted{11}{14}
\uarc45
\darc{15}{16}
\stlines{0/0,3/3,6/4,9/7,10/8,11/9,12/10,13/11,16/14}
\draw(-1,1)node[left]{$\lamb_i$};
\draw[|-|](-1,2)--(-1,0);
\node()at(0,2.5){\tiny$1$};
\node()at(4,2.5){\tiny$i$};
\node()at(9,2.5){\tiny$j$};
\node()at(16,2.5){\tiny$n$};
\end{scope}
\begin{scope}[shift={(0,-2)}]	
\uvs{0,3,4,7,8,9,10,11,14,15,16}
\lvs{0,3,4,7,8,9,12,13,14,15,16}
\udotted03
\udotted47
\udotted{11}{14}
\ddotted03
\ddotted47
\ddotted9{12}
\uarc9{10}
\darc{15}{16}
\stlines{0/0,3/3,4/4,7/7,8/8,11/9,14/12,15/13,16/14}
\draw(-1,1)node[left]{$\lamb_j$};
\draw[-|](-1,2)--(-1,0);
\end{scope}
\begin{scope}[shift={(0,-7)}]	
\uvs{0,3,4,5,6,9,10,11,12,13,16}
\lvs{0,3,4,5,6,9,10,11,14,15,16}
\udotted03
\udotted69
\udotted{13}{16}
\ddotted03
\ddotted69
\ddotted{11}{14}
\uarc{11}{12}
\darc{15}{16}
\stlines{0/0,3/3,4/4,5/5,6/6,9/9,10/10,13/11,16/14}
\draw(-1,1)node[left]{$\lamb_{j+2}$};
\draw[|-](-1,2)--(-1,0);
\node()at(0,2.5){\tiny$1$};
\node()at(4,2.5){\tiny$i$};
\node()at(9,2.5){\tiny$j$};
\node()at(16,2.5){\tiny$n$};
\end{scope}
\begin{scope}[shift={(0,-9)}]	
\uvs{0,3,4,5,6,9,10,11,14,15,16}
\lvs{0,3,4,7,8,9,12,13,14,15,16}
\udotted03
\udotted69
\udotted{11}{14}
\ddotted03
\ddotted47
\ddotted{9}{12}
\uarc45
\darc{15}{16}
\stlines{0/0,3/3,6/4,9/7,10/8,11/9,14/12,15/13,16/14}
\draw(-1,1)node[left]{$\lamb_i$};
\draw[|-|](-1,2)--(-1,0);
\end{scope}
\begin{scope}[shift={(21,0)}]	
\uvs{1,4,5,6,7,8,9,12,13}
\lvs{1,4,5,6,7,8,9,12,13}
\udotted14
\udotted9{12}
\ddotted14
\ddotted9{12}
\uarc67
\uarc89
\uarc{12}{13}
\darc67
\darc89
\darc{12}{13}
\stlines{1/1,4/4,5/5}
\draw(0,1)node[left]{$\lamb_k^i$};
\draw[|-|](0,2)--(0,0);
\node()at(1,2.5){\tiny$1$};
\node()at(6,2.5){\tiny$k$};
\node()at(13,2.5){\tiny$n$};
\end{scope}
\begin{scope}[shift={(21,-2)}]	
\uvs{1,4,5,6,7,8,9,12,13}
\lvs{1,4,5,6,7,10,11,12,13}
\udotted14
\udotted9{12}
\ddotted14
\ddotted7{10}
\uarc56
\darc{12}{13}
\stlines{1/1,4/4,7/5,8/6,9/7,12/10,13/11}
\draw(0,1)node[left]{$\lamb_{k-1}$};
\draw[-|](0,2)--(0,0);
\end{scope}
\begin{scope}[shift={(21,-7)}]	
\uvs{1,4,5,6,7,8,9,12,13}
\lvs{1,4,5,6,7,8,9,12,13}
\udotted14
\udotted9{12}
\ddotted14
\ddotted9{12}
\uarc67
\uarc89
\uarc{12}{13}
\darc67
\darc89
\darc{12}{13}
\stlines{1/1,4/4,5/5}
\draw(0,1)node[left]{$\lamb_k^i$};
\draw[|-|](0,2)--(0,0);
\node()at(1,2.5){\tiny$1$};
\node()at(6,2.5){\tiny$k$};
\node()at(13,2.5){\tiny$n$};
\end{scope}

\end{tikzpicture}
\caption{Left: relation \ref{L2}.  Right: relation \ref{L3}, writing $k=n-2i+1$.}
\label{f:L2}
\end{center}
\end{figure}

\pf
This amounts to checking that for each relation $(u,v)\in\Om_L$ we have $\ol u=\ol v$ in $\L_n$.  For~\ref{L1} and~\ref{L2}, this is easily accomplished diagrammatically; see Figure \ref{f:L2} for the latter.  For \ref{L3}, one may first show by a simple induction that $\bl(\lamb_{n-2i+1}^j) = (n-2i+1,n-2i+3,\ldots,n-2i+2j-1)$ for all $1\leq j\leq i$.  In particular, taking $j=i$, we have
\begin{equation}\label{e:llam}
\bl(\lamb_{n-2i+1}^i) = (n-2i+1,n-2i+3,\ldots,n-3,n-1).
\end{equation}
Since therefore $\rank(\lamb_{n-2i+1}^i)=n-2i$, we have $\br(\lamb_{n-2i+1}^i) = (n-2i+1,n-2i+3,\ldots,n-3,n-1)$ as well.  Thus,~\ref{L3} may be checked diagrammatically, as also shown in Figure \ref{f:L2}.
\epf

To complete the proof of Theorem \ref{t:LR}, we must prove the reverse inclusion: $\ker(\phi_L)\sub{\sim_L}$.  First we need a simple consequence of the relations.  Note that \ref{L1} says $\lam_{n-1}$ is a right identity for each $\lam_j$.  Since the $i=1$ case of \ref{L3} says $\lam_{n-1}\lam_{n-2}=\lam_{n-1}$, it follows that $\lam_{n-2}$ is also a right identity:

\begin{lemma}\label{l:n-2}
If $1\leq i\leq n-1$ and $n-2\leq j\leq n-1$, then $\lam_i\lam_j\sim_L\lam_i$.  \epfres
\end{lemma}

For a tuple $\bx=(x_1,\ldots,x_k)\in T_n$ (as defined in Section \ref{ss:tuples}), we define the words
\[
\lam_\bx = \lam_{x_1}\cdots\lam_{x_k} \in L^* \AND \rho_\bx = \rho_{x_k}\cdots\rho_{x_1}\in R^*.
\]
Note that when $\bx=\emptyset$ is the empty tuple, $\lam_\emptyset=\rho_\emptyset=1$ is the empty word.  The next lemma refers to the elements $\lamb_\bx=\lamb_{x_1}\cdots\lamb_{x_k}\in\L_n$.

\begin{lemma}\label{l:lambbx}
If $\bx\in T_n$, then $\bl(\lamb_\bx)=\bx$.
\end{lemma}

\pf
This is clear if $|\bx|\leq1$.  Otherwise, write $\bx=(x_1,\ldots,x_k)$, where $k\geq2$.  Put $\by=(x_1,\ldots,x_{k-1})$, noting that $\by\in T_n$, and that $\bl(\lamb_\by)=\by$ by induction.  Since $\lamb_\by\in\L_n$ has $k-1$ non-transversals, we have $\codom(\lamb_\by)=\{1,\ldots,n-2k+2\}$.  Since $x_k\leq n-2k+1$, it follows that $x_k$ and $x_k+1$ both belong to $\codom(\lamb_\by)$.  
Since $x_k<x_{k-1}<\cdots<x_1$, each of $\lamb_{x_1},\ldots,\lamb_{x_{k-1}}$ contains a transversal joining $P_{x_k}$ to~$P_{x_k}'$; so too therefore does $\lamb_{x_1}\cdots\lamb_{x_{k-1}}=\lamb_\by$.  It follows that $\lamb_\bx=\lamb_\by\lamb_{x_k}$ contains a non-transversal joining $P_{x_k}$ to~$P_y$ for some $y>x_k$.  Together with the fact that each non-transversal of $\lamb_\by$ is contained in $\lamb_\by\lamb_{x_k}=\lamb_\bx$, it follows that~$\bl(\lamb_\bx)=(x_1,\ldots,x_{k-1},x_k)=\bx$.
\epf

The following important consequence will be used in the next section:

\begin{cor}\label{c:*}
If $\bx,\by\in T_n$ and $|\bx|=|\by|$, then $\bl(\lamb_\bx\rhob_\by)=\bx$ and $\br(\lamb_\bx\rhob_\by)=\by$.
\end{cor}

\pf
Write $k=|\bx|=|\by|$.  Then $\codom(\lamb_\bx)=\{1,\ldots,n-2k\}=\dom(\rhob_\by)$, so ${\rank(\lamb_\bx\rhob_\by)=n-2k}$.  Thus, $\lamb_\bx\rhob_\by$ has precisely the same upper (respectively, lower) non-transversals as $\lamb_\bx$ (respectively, $\rhob_\by$), and so $\bl(\lamb_\bx\rhob_\by)=\bl(\lamb_\bx)$ and ${\br(\lamb_\bx\rhob_\by)=\br(\rhob_\by)}$.  The result now follows from Lemma \ref{l:lambbx} and its dual.
\epf

We now wish to show that every word over $L$ is $\sim_L$-equivalent to one of the form $\lam_\bx$ with $\bx\in T_n$.  This is accomplished in Lemma \ref{l:lambx} below, for whose inductive proof we require some technical lemmas analysing products of the form $\lam_\bx\lam_j$, where $\bx\in T_n$ and $1\leq j\leq n-1$.

\begin{lemma}\label{l:xj1}
If $\bx\in T_n$ with $k=|\bx|\geq1$, then $\lam_\bx\lam_j\sim_L\lam_\bx$ for any $j\geq n-2k$.
\end{lemma}

\pf
If $j\geq n-2$, then the result follows from Lemma \ref{l:n-2}.  This includes the $k=1$ case, so we now assume that $k\geq2$ and $j\leq n-3$, and we proceed by induction on $k$.  Write $\bx=(x_1,\ldots,x_k)$, and also let $\by=(x_1,\ldots,x_{k-1})$.  If $j\geq n-2k+1$, then since $x_k\leq n-2k+1\leq j\leq n-3$, we may apply~\ref{L2} and induction to obtain $\lam_\bx\lam_j = \lam_\by\lam_{x_k}\lam_j \sim_L \lam_\by\lam_{j+2}\lam_{x_k} \sim_L \lam_\by\lam_{x_k} =\lam_\bx$.
This leaves only the case in which $j=n-2k$.  For this we use the $j=n-2k+1$ case (just proved) and \ref{L3} to obtain $\lam_\bx\lam_{n-2k} \sim_L \lam_\bx\lam_{n-2k+1}^k\lam_{n-2k} \sim_L \lam_\bx\lam_{n-2k+1}^k \sim_L \lam_\bx$
\epf

For the next statement, we define $\min(\bx)$ to be the minimum entry of $\bx\in T_n$ if $|\bx|\geq1$, and also $\min(\emptyset)=\infty$.

\begin{lemma}\label{l:xj3}
If $\bx\in T_n$ with $k=|\bx|\geq0$, and if $j\leq n-2k-1$, then $\lam_\bx\lam_j\sim_L\lam_\by$ for some $\by\in T_n$ with $|\by|=k+1$ and $\min(\by) = \min(\min(\bx),j)$.
\end{lemma}

\pf
We use induction on $k$.  If $k=0$, then $\lam_\bx\lam_j=\lam_j=\lam_\by$, where $\by=(j)\in T_n$.  Now suppose $k\geq1$ and write $\bx=(x_1,\ldots,x_k)$, noting that $\min(\bx)=x_k$.
If $x_k>j$, then $\lam_\bx\lam_j=\lam_\by$ where $\by=(x_1,\ldots,x_k,j)\in T_n$.  So for the rest of the proof we will assume that $x_k\leq j$, and we define $\bu=(x_1,\ldots,x_{k-1})$, which is empty if $k=1$.  Since $x_k\leq j\leq n-2k-1\leq n-3$ (as $k\geq1$), we apply~\ref{L2} and then induction to obtain $\lam_\bx\lam_j = \lam_\bu \lam_{x_k}\lam_j \sim_L \lam_\bu \lam_{j+2} \lam_{x_k} \sim_L \lam_\bv \lam_{x_k}$ for some $\bv\in T_n$ with $|\bv|=k$ and $\min(\bv)=\min(\min(\bu),j+2)$.

Write $\bv=(v_1,\ldots,v_k)$, noting that $v_k=\min(\bv)$ is either $j+2$ or $\min(\bu)$; also note that $\min(\bu)=\infty$ if $k=1$, or else $\min(\bu)=x_{k-1}$.  In any case, we have $v_k>x_k$ (since $j+2\geq x_k+2$, and since $x_{k-1}>x_k$ if $k\geq2$), and so $\lam_\bv\lam_{x_k}=\lam_\by$, where $\by=(v_1,\ldots,v_k,x_k)\in T_n$.  (Note that $x_k\leq j\leq n-2k-1$.)
\epf

Here is the most important technical lemma of this section.  

\begin{lemma}\label{l:lambx}
If $w\in L^*$, then $w\sim_L \lam_\bx$, where $\bx=\bl(\ol w)$.
\end{lemma}

\pf
First note that it suffices to show that $w\sim_L\lam_\bx$ for some $\bx\in T_n$, since we will then have $\bl(\ol w)=\bl(\lamb_\bx)=\bx$, by Lemma \ref{l:lambbx}.  To do this, we use induction on $l$, the length of $w$.  If $l\leq 1$, then~$w$ is already of the form $\lam_\bx$, where $|\bx|\leq1$.  If $l\geq2$, then $w=u\lam_j$ for some $u\in L^*$ of length $l-1$, and some $1\leq j\leq n-1$.  By induction, we have $u\sim_L\lam_\by$ for some $\by\in T_n$ with $k=|\by|\geq1$.  If $j\geq n-2k$, then Lemma \ref{l:xj1} gives $w=u\lam_j\sim_L\lam_\by\lam_j\sim_L\lam_\by$.  If $j\leq n-2k-1$, then Lemma \ref{l:xj3} gives $w\sim_L\lam_\by\lam_j\sim_L\lam_\bx$ for some $\bx\in T_n$.
\epf

We may now tie together the loose ends.

\pf[\bf Proof of Theorem \ref{t:LR}]
Recall that we only need to give the proof for $\L_n$.  For this, it remains to show that $\ker(\phi_L)\sub{\sim_L}$.  To do so, let $(u,v)\in\ker(\phi_L)$; so $u,v\in L^*$ and $\ol u=\ol v$.  By Lemma \ref{l:lambx}, we have $u\sim_L\lam_\bx$, where $\bx=\bl(\ol u)$.  Since $\ol v=\ol u$, we also have $v\sim_L\lam_\bx$, and so $u\sim_L\lam_\bx\sim_L v$.
\epf

\begin{rem}\label{r:nf}
The proofs of Lemmas \ref{l:xj1}--\ref{l:lambx} give an algorithm for transforming an arbitrary word over $L$ into a normal form $\lam_\bx$, for $\bx\in T_n$.
\end{rem}

\section{First presentation for $\TL_n$}\label{s:TL1}

The previous section gave presentations for the monoids $\L_n$ and $\R_n$ (see Theorem \ref{t:LR}), and we now use these to complete the proof of Theorem \ref{t:TL1}.

\begin{prop}\label{p:TL_LR}
If $\al\in\TL_n$, then $\al=\lamb_\bx\rhob_\by$, where $\bx=\bl(\al)$ and $\by=\br(\al)$.  Consequently,
\[
\TL_n=\L_n\R_n=\la\ol L\cup\ol R\ra.
\]
\end{prop}

\begin{figure}[h]
\begin{center}
\begin{tikzpicture}[xscale=.45,yscale=0.45]
\ublueboxes{1/3,5/6,8/10,15/16,18/20}
\lblueboxes{1/2,4/8,10/11,14/15,17/20}
\uvs{1,3,4,5,6,7,8,10,15,16,17,18,20}
\lvs{1,2,3,4,8,9,10,11,14,15,16,17,20}
\stlines{4/3,7/9,17/16}
\udottedsms{10/15}
\ldottedsms{11/14}
  \draw[|-|] (.2,0)--(.2,2);
  \draw(.2,1)node[left]{$\al$};
\begin{scope}[shift={(0,-8)}]	
\uublueboxes{1/3,5/6,8/10,15/16,18/20}
\uuvs{1,3,4,5,6,7,8,10,15,16,17,18,20}
\lvs{1,2,7,8,9,10,11,19,20}
\darcx89{.3}
\darcx{10}{11}{.3}
\darcx{19}{20}{.3}
\ustlines{4/1,7/2,17/7}
\uudottedsms{10/15}
\ldottedsms{2/7,11/19}
  \draw[-|] (.2,0)--(.2,4);
  \draw(.2,2)node[left]{$\be$};
\end{scope}
\begin{scope}[shift={(0,-12)}]	
\lblueboxes{1/2,4/8,10/11,14/15,17/20}
\uuvs{1,2,7,8,9,10,11,19,20}
\lvs{1,2,3,4,8,9,10,11,14,15,16,17,20}
\uuarcx89{.3}
\uuarcx{10}{11}{.3}
\uuarcx{19}{20}{.3}
\ustlines{1/3,2/9,7/16}
\uudottedsms{2/7,11/19}
\ldottedsms{11/14}
  \draw[|-|] (.2,0)--(.2,4);
  \draw(.2,2)node[left]{$\ga$};
\end{scope}
\end{tikzpicture}
\caption{Partitions constructed during the proof of Proposition \ref{p:TL_LR}; note that $\al=\be\ga$.}
\label{f:TL_LR}
\end{center}
\end{figure}
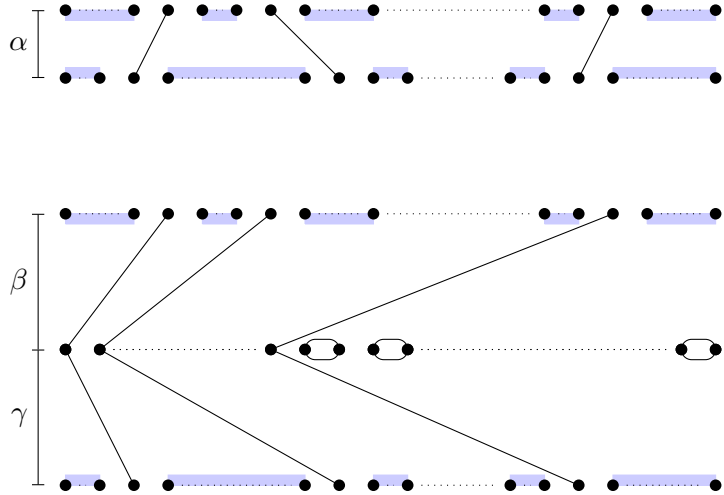

\pf
In light of Proposition \ref{p:L}, it suffices to prove the first statement.  Clearly we have $\al=\be\ga$, where $\be\in\L_n$ and $\ga\in\R_n$ are such that $\be$ (respectively, $\ga$) has the same upper (respectively, lower) non-transversals as $\al$; see Figure \ref{f:TL_LR}.  Since $\be$ has the same upper non-transversals as $\al$, we have $\bl(\be)=\bl(\al)=\bx=\bl(\lamb_\bx)$, where we used Lemma \ref{l:lambbx} in the last step.  It follows from Corollary \ref{c:LR}\ref{c:LR1} that $\be=\lamb_\bx$.  A symmetrical argument shows that $\ga=\rhob_\by$.  Thus, $\al=\be\ga=\lamb_\bx\rhob_\by$.
\epf

The next simple consequence refers to the partitions $\bp(\al)$ as defined in Section \ref{ss:tuples}.  The equivalence of \ref{LR1} and \ref{LR2} was stated without proof in \cite[Remark 3]{BDP2002}.

\begin{cor}\label{c:TL_LR}
For $\al,\be\in\TL_n$, the following are equivalent:
\ben\bmc3
\item \label{LR1} $\al=\be$,
\item \label{LR2} $\bp(\al)=\bp(\be)$,
\item \label{LR3} $\bl(\al)=\bl(\be)$ and $\br(\al)=\br(\be)$.
\emc\een
\end{cor}

\pf
The implications \ref{LR1}$\implies$\ref{LR2} and \ref{LR2}$\implies$\ref{LR3} are clear, while \ref{LR3}$\implies$\ref{LR1} follows from Proposition~\ref{p:TL_LR}: if $\bl(\al)=\bl(\be)=\bx$ and $\br(\al)=\br(\be)=\by$, then $\al=\lamb_\bx\rhob_\by=\be$.
\epf

It also follows from Proposition \ref{p:TL_LR} that the morphism $\phi:(L\cup R)^*\to\TL_n$ defined in Section \ref{ss:statements} is surjective.  To complete the proof of Theorem \ref{t:TL1}, we must show that $\ker(\phi)=\Om^\sharp$, where $\Om$ is the set of relations \ref{L1}--\ref{L3}, \ref{R1}--\ref{R3} and \ref{RL1}--\ref{RL3}.  Throughout this section, we write ${\sim}=\Om^\sharp$.  

\begin{lemma}
We have ${\sim}\sub\ker(\phi)$.
\end{lemma}

\pf
This again amounts to a simple diagrammatic check; see Figure \ref{f:RL3} for \ref{RL3}.
\epf

\begin{figure}[ht]
\begin{center}
\begin{tikzpicture}[scale=.42]
\begin{scope}[shift={(0,0)}]	
\uvs{1,4,5,8,9,10,11,14,15,16}
\lvs{1,4,5,6,7,10,11,12,13,16}
\udotted14
\udotted58
\udotted{11}{14}
\ddotted14
\ddotted7{10}
\ddotted{13}{16}
\uarc{15}{16}
\darc56
\stlines{1/1,4/4,5/7,8/10,9/11,10/12,11/13,14/16}
\draw(0,1)node[left]{$\rhob_i$};
\draw[|-|](0,2)--(0,0);
\node()at(1,2.5){\tiny$1$};
\node()at(5,2.5){\tiny$i$};
\node()at(11,2.5){\tiny$j$};
\node()at(16,2.5){\tiny$n$};
\end{scope}
\begin{scope}[shift={(0,-2)}]	
\uvs{1,4,5,6,7,10,11,12,13,16}
\lvs{1,4,5,6,7,10,11,14,15,16}
\udotted14
\udotted7{10}
\udotted{13}{16}
\ddotted14
\ddotted7{10}
\ddotted{11}{14}
\uarc{11}{12}
\darc{15}{16}
\stlines{1/1,4/4,5/5,6/6,7/7,10/10,13/11,16/14}
\draw(0,1)node[left]{$\lamb_j$};
\draw[-|](0,2)--(0,0);
\end{scope}
\begin{scope}[shift={(22,0)}]	
\uvs{1,4,5,8,9,10,11,14,15,16}
\lvs{1,4,5,8,9,10,11,14,15,16}
\udotted14
\udotted58
\udotted{11}{14}
\ddotted14
\ddotted58
\ddotted{11}{14}
\uarc{15}{16}
\darc{15}{16}
\stlines{1/1,4/4,5/5,8/8,9/9,10/10,11/11,14/14}
\draw(0,1)node[left]{$\lamb_{n-1}$};
\draw[|-|](0,2)--(0,0);
\node()at(1,2.5){\tiny$1$};
\node()at(5,2.5){\tiny$i$};
\node()at(11,2.5){\tiny$j$};
\node()at(16,2.5){\tiny$n$};
\end{scope}
\begin{scope}[shift={(22,-2)}]	
\uvs{1,4,5,8,9,10,11,14,15,16}
\lvs{1,4,5,8,9,12,13,14,15,16}
\udotted14
\udotted58
\udotted{11}{14}
\ddotted14
\ddotted58
\ddotted9{12}
\uarc9{10}
\darc{15}{16}
\stlines{1/1,4/4,5/5,8/8,11/9,14/12,15/13,16/14}
\draw(0,1)node[left]{$\lamb_{j-2}$};
\draw[-|](0,2)--(0,0);
\end{scope}
\begin{scope}[shift={(22,-4)}]	
\uvs{1,4,5,8,9,12,13,14,15,16}
\lvs{1,4,5,6,7,10,11,14,15,16}
\udotted14
\udotted58
\udotted9{12}
\ddotted14
\ddotted7{10}
\ddotted{11}{14}
\uarc{15}{16}
\darc56
\stlines{1/1,4/4,5/7,8/10,9/11,12/14,13/15,14/16}
\draw(0,1)node[left]{$\rhob_i$};
\draw[-|](0,2)--(0,0);
\end{scope}
\end{tikzpicture}
\caption{Relation \ref{RL3}.}
\label{f:RL3}
\end{center}
\end{figure}

We now work towards the reverse inclusion, beginning with an obvious consequence of \ref{RL1}--\ref{RL3}:

\begin{lemma}\label{l:uv}
For any $w\in(L\cup R)^*$, we have $w\sim uv$ for some $u\in L^*$ and $v\in R^*$.  \epfres
\end{lemma}

It follows from Lemma \ref{l:uv}, along with Lemma \ref{l:lambx} and its dual, that any word over $L\cup R$ is $\sim$-equivalent to $\lam_\bx\rho_\by$ for some $\bx,\by\in T_n$.  The main step remaining is to show that such a factorisation exists with $|\bx|=|\by|$, and this is achieved in Lemma \ref{l:wxy}.  Roughly speaking, we use the longer of~$\lam_\bx$ or~$\rho_\by$ to lengthen the shorter, and the next lemma provides the precise mechanism to do this.

\begin{lemma}\label{l:lrlr}
If $\bx\in T_n$ with $k=|\bx|\geq1$, then $\lam_\bx\sim\lam_\bx\rho_{n-2k+1}$ and $\rho_\bx\sim\lam_{n-2k+1}\rho_\bx$.
\end{lemma}

\pf
By the asymmetry in relations \ref{RL1}--\ref{RL3}, the two statements are not dual, so both must be treated.  We proceed by induction on $k$.  Write $\bx=(x_1,\ldots,x_k)$.

If $k=1$, then $\lam_\bx\sim\lam_\bx\lam_{n-1}\sim\lam_\bx\rho_{n-1}$ and $\rho_\bx\sim\rho_{n-1}\rho_\bx\sim\lam_{n-1}\rho_\bx$, using \ref{L1}, \ref{R1} and \ref{RL2}.  Now suppose $k\geq2$, and write $\by=(x_1,\ldots,x_{k-1})\in T_n$.  Then
\[
\lam_\bx = \lam_\by\lam_{x_k} \sim \lam_\by\rho_{n-2k+3}\lam_{x_k} 
\sim \lam_\by\lam_{n-1}\lam_{x_k}\rho_{n-2k+1} 
\sim \lam_\by\lam_{x_k}\rho_{n-2k+1} 
= \lam_\bx\rho_{n-2k+1},
\]
where we first used induction, and then \ref{RL1} and \ref{L1}, noting that $x_k\leq n-2k+1$.  We also have
\begin{align*}
\rho_\bx = \rho_{x_k} \rho_\by &\sim \rho_{x_k} \lam_{n-2k+3} \rho_\by &&\text{by induction}\\
&\sim \lam_{n-1} \lam_{n-2k+1} \rho_{x_k}\rho_\by &&\text{by \ref{RL3}}\\
&\sim  \lam_{n-2k+1}\lam_{n-3} \rho_{x_k}\rho_\by &&\text{by \ref{L2}, with $i=n-2k+1$ and $j=n-3$}\\
&\sim  \lam_{n-2k+1}\lam_{n-1}\lam_{n-3} \rho_{x_k}\rho_\by &&\text{by \ref{L1}}\\
&\sim  \lam_{n-2k+1} \rho_{x_k}\lam_{n-1}\rho_\by &&\text{by \ref{RL3}, with $i=x_k$ and $j=n-1$}\\
&\sim  \lam_{n-2k+1} \rho_{x_k}\rho_{n-1}\rho_\by &&\text{by \ref{RL2}}\\
&\sim  \lam_{n-2k+1} \rho_{x_k}\rho_\by=\lam_{n-2k+1}\rho_\bx &&\text{by \ref{R1}.} \qedhere
\end{align*}
\epf

Here are the promised normal forms:

\begin{lemma}\label{l:wxy}
If $w\in(L\cup R)^*$, then $w\sim\lam_\bx\rho_\by$, where $\bx=\bl(\ol w)$ and $\by=\br(\ol w)$.
\end{lemma}

\pf
It is enough to show that $w\sim\lam_\bx\rho_\by$ for some $\bx,\by\in T_n$ with $|\bx|=|\by|$.  Indeed, then we will have $\bl(\ol w)=\bl(\lamb_\bx\rhob_\by)=\bx$, by Corollary \ref{c:*}, and similarly $\br(\ol w)=\by$.

If $w=1$ then $w=\lam_\emptyset\rho_\emptyset$, so we assume that $w\not=1$.  By Lemma \ref{l:uv} we have $w\sim uv$ for some $u\in L^*$ and $v\in R^*$.  Write $k=|\bl(\ol u)|$ and $l=|\br(\ol v)|$, and assume that $k\geq l$, the other case being symmetrical.  Since $w\not=1$, we have $k\geq1$.  Since $\ol u\in\L_n$ and $\ol v\in\R_n$ have $k$ and~$l$ non-transversals, respectively, we have $\codom(\ol u)=\{1,\ldots,n-2k\}\sub\{1,\ldots,n-2l\}=\dom(\ol v)$.  It follows that
\begin{equation}\label{e:rankw}
\rank(\ol w)=\rank(\ol u \;\!\;\! \ol v)=\rank(\ol u)=n-2k.
\end{equation}
Writing $\bx=\bl(\ol u)$ and $\bz=\br(\ol v)$, Lemma \ref{l:lambx} and its dual give $u\sim\lam_\bx$ and $v\sim\rho_\bz$.  Together with~$k$ applications of Lemma \ref{l:lrlr}, it follows that $w \sim uv \sim \lam_\bx\rho_\bz \sim \lam_\bx\rho_{n-2k+1}^k\rho_\bz$.  Combining this with \eqref{e:llam} and~\eqref{e:rankw}, it follows that
\[
\rank(\rhob_{n-2k+1}^k\rhob_\bz) \leq \rank(\rhob_{n-2k+1}^k) = n-2k = \rank(\ol w) = \rank(\lamb_\bx\rhob_{n-2k+1}^k\rhob_\bz) \leq \rank(\rhob_{n-2k+1}^k\rhob_\bz) ,
\]
so that in fact $\rank(\rhob_{n-2k+1}^k\rhob_\bz)=n-2k$.  Thus,~$\rhob_{n-2k+1}^k\rhob_\bz$ has $k$ non-transversals and so, writing $\by=\br(\rhob_{n-2k+1}^k\rhob_\bz)$, we have $|\by|=k=|\bx|$, and the dual of Lemma \ref{l:lambx} gives $\rho_{n-2k+1}^k\rho_\bz\sim\rho_\by$.  Putting everything together, we have $w \sim \lam_\bx\rho_{n-2k+1}^k\rho_\bz \sim \lam_\bx\rho_\by$.
\epf

As in Remark \ref{r:nf}, the above lemmas lead to an algorithm to compute the normal form of an arbitrary word over $L\cup R$.  
We may now complete the proof of our main new result.

\pf[\bf Proof of Theorem \ref{t:TL1}]
It remains to show that $\ker(\phi)\sub{\sim}$.  To do so, let $(u,v)\in\ker(\phi)$; so $u,v\in(L\cup R)^*$ and $\ol u=\ol v$.  By Lemma \ref{l:wxy}, we have $u\sim\lam_\bx\rho_\by$, where $\bx=\bl(\ol u)$ and $\by=\br(\ol u)$.  Since $\ol v=\ol u$, we also have $v\sim\lam_\bx\rho_\by$, and so $u\sim\lam_\bx\rho_\by\sim v$.
\epf

\section{Second presentation for $\TL_n$}\label{s:TL2}

To prove Theorem \ref{t:TL2}, we need to show that $\psi:E^*\to\TL_n:e_i\mt\ol e_i$ is surjective, and that $\ker(\psi)=\Xi^\sharp$.  

\begin{prop}\label{p:TL}
We have $\TL_n=\la\ol E\ra$.  Consequently, $\psi$ is surjective.
\end{prop}

\pf
Since $\lamb_i = \ol e_i\ol e_{i+1}\cdots\ol e_{n-1}$ and $\rhob_i = \ol e_{n-1}\cdots\ol e_{i+1}\ol e_i$ for all $i$, this follows from Proposition~\ref{p:TL_LR}.
\epf

For the rest of this section we write ${\approx}=\Xi^\sharp$.  The next result is proved in the usual way.  

\begin{lemma}\label{l:Xi}
We have ${\approx}\sub\ker(\psi)$.  \epfres
\end{lemma}

As ever, proving the reverse inclusion is more of a challenge, though we will be greatly aided by the fact that we already have the presentation $\pres{L\cup R}\Om$ from Theorem \ref{t:TL1}.

For $1\leq i\leq n-1$, we define the following words from $E^*$:
\[
\wh\lam_i =  e_i e_{i+1}\cdots e_{n-1} \AND \wh\rho_i =  e_{n-1}\cdots e_{i+1} e_i.
\]
We extend this to a morphism $(L\cup R)^*\to E^*:w\mt \wh w$, and we note that $\ol{\ \wh w\ } = \ol w$ for all $w\in (L\cup R)^*$, meaning that $\wh w\psi = w\phi$ for all such $w$.

\begin{lemma}\label{l:wh}
For any $1\leq i\leq n-1$, we have $\wh\lam_i\wh\rho_i\approx e_i$. 
\end{lemma}

\pf
We use reverse induction on $i$.  For $i=n-1$, we have $\wh\lam_{n-1}\wh\rho_{n-1}=e_{n-1}e_{n-1}\approx e_{n-1}$, by~\ref{E1}.  If $i\leq n-2$, then by induction and \ref{E3}, we have $\wh\lam_i\wh\rho_i = e_i\wh\lam_{i+1}\wh\rho_{i+1}e_i \approx e_i e_{i+1}e_i \approx e_i$.
\epf

The next result follows immediately:

\begin{cor}\label{c:wh}
For any $w\in E^*$, we have $w\approx\wh u$ for some $u\in(L\cup R)^*$.  \epfres
\end{cor}

The main technical result we need provides the link between the relations in the two presentations:

\begin{lemma}\label{l:Om}
For any relation $(u,v)$ from $\Om$, we have $\wh u\approx\wh v$.
\end{lemma}

\pf
For \ref{L1} this follows immediately from \ref{E1}.
For \ref{L2}, if $i\leq j\leq n-3$, then
\begin{align*}
\wh\lam_{j+2}\wh\lam_i &= e_{j+2}\cdots e_{n-1} \cdot e_i\cdots e_{j-1}\cdot e_j\cdot e_{j+1}\cdots e_{n-1} \\
&\approx e_i\cdots e_{j-1}\cdot e_j\cdot e_{j+2}\cdots e_{n-1} \cdot  e_{j+1}\cdots e_{n-1} &&\text{by \ref{E2}}\\
&\approx e_i\cdots e_{j-1}\cdot e_je_{j+1}e_j\cdot e_{j+2}\cdots e_{n-1} \cdot  e_{j+1}\cdots e_{n-1} &&\text{by \ref{E3}}\\
&\approx e_i\cdots e_{j-1}\cdot e_je_{j+1}\cdot e_{j+2}\cdots e_{n-1} \cdot  e_je_{j+1}\cdots e_{n-1} = \wh\lam_i\wh\lam_j &&\text{by \ref{E2}.}
\end{align*}
We treat \ref{L3} by induction.  For $i=1$, \ref{E3} gives $\wh\lam_{n-1}\wh\lam_{n-2} = e_{n-1}e_{n-2}e_{n-1} \approx e_{n-1} = \wh\lam_{n-1}$.  For~$i\geq2$,
\begin{align*}
\wh\lam_{n-2i+1}^i\wh\lam_{n-2i} &= \wh\lam_{n-2i+1}\cdot\wh\lam_{n-2i+1}^{i-1}\cdot\wh\lam_{n-2i} \\
&\approx \wh\lam_{n-2i+3}^{i-1}\cdot\wh\lam_{n-2i+1}\cdot\wh\lam_{n-2i} &&\text{by \ref{L2}, proved above}\\
&= \wh\lam_{n-2i+3}^{i-1}\cdot e_{n-2i+1}\wh\lam_{n-2i+2}\cdot\wh\lam_{n-2i} \\
&\approx e_{n-2i+1}\cdot \wh\lam_{n-2i+3}^{i-1} \wh\lam_{n-2i+2}\cdot\wh\lam_{n-2i} &&\text{by \ref{E2}}\\
&\approx e_{n-2i+1}\cdot \wh\lam_{n-2i+3}^{i-1} \cdot\wh\lam_{n-2i} &&\text{by induction}\\
&\approx  \wh\lam_{n-2i+3}^{i-1} \cdot e_{n-2i+1}\wh\lam_{n-2i} &&\text{by \ref{E2}}\\
&= \wh\lam_{n-2i+3}^{i-1} \cdot e_{n-2i+1}e_{n-2i}e_{n-2i+1}\wh\lam_{n-2i+2} \\
&\approx \wh\lam_{n-2i+3}^{i-1} \cdot e_{n-2i+1}\wh\lam_{n-2i+2} &&\text{by \ref{E3}}\\
&= \wh\lam_{n-2i+3}^{i-1} \cdot \wh\lam_{n-2i+1} \\
&\approx \wh\lam_{n-2i+1} \cdot \wh\lam_{n-2i+1}^{i-1} = \wh\lam_{n-2i+1}^i &&\text{by \ref{L2} again.}
\end{align*}
We have dealt with relations \ref{L1}--\ref{L3}.  Relations \ref{R1}--\ref{R3} are dual.

For \ref{RL2}, first note that $\wh\rho_{n-1}\wh\lam_{n-1}=e_{n-1}e_{n-1}\approx e_{n-1}=\wh\lam_{n-1}=\wh\rho_{n-1}$, while if $i\leq n-2$, then
\begin{align*}
\wh\rho_i\wh\lam_i &= e_{n-1}\cdots e_{i+2}e_{i+1}e_i\cdot e_ie_{i+1}e_{i+2}\cdots e_{n-1} \\
&\approx e_{n-1}\cdots e_{i+2}e_{i+1}e_i\cdot e_{i+1}e_{i+2}\cdots e_{n-1} &&\text{by \ref{E1}}\\
&\approx e_{n-1}\cdots e_{i+2}\cdot e_{i+1}e_{i+2}\cdots e_{n-1} &&\text{by \ref{E3}}\\
&\approx e_{n-1}\cdots e_{i+2}e_{i+1}\cdot e_{i+1}e_{i+2}\cdots e_{n-1} &&\text{by \ref{E1}}\\
&= \wh\rho_{i+1}\wh\lam_{i+1} \approx \wh\lam_{n-1} &&\text{by induction.}
\end{align*}
The second line in the previous calculation also shows that $\wh\rho_i\wh\lam_i\approx\wh\rho_i\wh\lam_{i+1}=\wh\rho_{i+1}\wh\lam_i$ for $i\leq n-2$, and this finishes the proof for \ref{RL2}.
For \ref{RL1}, if $j\leq i-2$, then
\begin{align*}
\wh\rho_i\wh\lam_j &= e_{n-1}\cdots e_i\cdot e_j\cdots e_{i-3}e_{i-2}\cdot e_{i-1}e_i\cdots e_{n-1} \\
&\approx e_j\cdots e_{i-3}e_{i-2}\cdot e_{n-1}\cdots e_i\cdot e_{i-1}e_i\cdots e_{n-1} &&\text{by \ref{E2}}\\
&= e_j\cdots e_{i-3}e_{i-2}\cdot \wh\rho_i\wh\lam_{i-1} \\
&\approx e_j\cdots e_{i-3}e_{i-2}\cdot \wh\lam_{n-1} &&\text{by \ref{RL2}, proved above}\\
&= e_j\cdots e_{i-3}e_{i-2}\cdot e_{n-1} \\
&\approx e_{n-1}\cdot e_j\cdots e_{i-3}e_{i-2}  &&\text{by \ref{E2}} \\
&\approx e_{n-1}\cdot e_j\cdots e_{i-3}\wh\lam_{i-2}\wh\rho_{i-2} = \wh\lam_{n-1} \wh\lam_j\wh\rho_{i-2} &&\text{by Lemma \ref{l:wh}.}
\end{align*}
Relation \ref{RL3} is treated in almost identical fashion to \ref{RL1}.
\epf

We now have everything we need to complete the proof of the second main result.

\pf[{\bf Proof of Theorem \ref{t:TL2}}]
It remains to show that $\ker(\psi)\sub{\approx}$.  To do so, let $(w_1,w_2)\in\ker(\psi)$.  So $w_1,w_2\in E^*$ and $\ol w_1=\ol w_2$.  By Corollary \ref{c:wh}, we have $w_1\approx\wh u$ and $w_2\approx \wh v$ for some $u,v\in(L\cup R)^*$, and we note that $\ol u = \ol{\ \wh u\ } = \ol w_1 = \ol w_2 = \ol{\ \wh v\ } = \ol v$.  
Since $\ol u=\ol v$, it follows from Theorem \ref{t:TL1} that there is a sequence $u = u_0 \to u_1 \to \cdots \to u_k = v$, where for each $1\leq i\leq k$, $u_i\in(L\cup R)^*$ is obtained from $u_{i-1}$ by a single application of one of the relations from $\Om$.  It then follows from Lemma \ref{l:Om} that $\wh u = \wh u_0 \approx \wh u_1 \approx\cdots\approx\wh u_k = \wh v$.
But then $w_1\approx\wh u\approx\wh v\approx w_2$.
\epf

\footnotesize
\def\bibspacing{-1.1pt}
\bibliography{biblio}

\def\cprime{$'$}
\begin{thebibliography}{10}

\bibitem{Abramsky2008}
S.~Abramsky.
\newblock Temperley-{L}ieb algebra: from knot theory to logic and computation
  via quantum mechanics.
\newblock In {\em Mathematics of quantum computation and quantum technology},
  Chapman \& Hall/CRC Appl. Math. Nonlinear Sci. Ser., pages 515--558. Chapman
  \& Hall/CRC, Boca Raton, FL, 2008.

\bibitem{Auinger2014}
K.~Auinger.
\newblock Pseudovarieties generated by {B}rauer type monoids.
\newblock {\em Forum Math.}, 26(1):1--24, 2014.

\bibitem{BDP2002}
M.~Borisavljevi{\'c}, K.~Do{\v{s}}en, and Z.~Petri{\'c}.
\newblock Kauffman monoids.
\newblock {\em J. Knot Theory Ramifications}, 11(2):127--143, 2002.

\bibitem{Brauer1937}
R.~Brauer.
\newblock On algebras which are connected with the semisimple continuous
  groups.
\newblock {\em Ann. of Math. (2)}, 38(4):857--872, 1937.

\bibitem{DEEFHHLM2019}
I.~Dolinka, J.~East, A.~Evangelou, D.~FitzGerald, N.~Ham, J.~Hyde, N.~Loughlin,
  and J.~D. Mitchell.
\newblock Enumeration of idempotents in planar diagram monoids.
\newblock {\em J. Algebra}, 522:351--385, 2019.

\bibitem{JEinsn}
J.~East.
\newblock A presentation of the singular part of the symmetric inverse monoid.
\newblock {\em Comm. Algebra}, 34(5):1671--1689, 2006.

\bibitem{JEgrpm}
J.~East.
\newblock Generators and relations for partition monoids and algebras.
\newblock {\em J. Algebra}, 339:1--26, 2011.

\bibitem{EG2017}
J.~East and R.~D. Gray.
\newblock Diagram monoids and {G}raham--{H}oughton graphs: {I}dempotents and
  generating sets of ideals.
\newblock {\em J. Combin. Theory Ser. A}, 146:63--128, 2017.

\bibitem{EMRT2018}
J.~East, J.~D. Mitchell, N.~Ru\v{s}kuc, and M.~Torpey.
\newblock Congruence lattices of finite diagram monoids.
\newblock {\em Adv. Math.}, 333:931--1003, 2018.

\bibitem{Jones1983_2}
V.~F.~R. Jones.
\newblock Index for subfactors.
\newblock {\em Invent. Math.}, 72(1):1--25, 1983.

\bibitem{Jones1987}
V.~F.~R. Jones.
\newblock Hecke algebra representations of braid groups and link polynomials.
\newblock {\em Ann. of Math. (2)}, 126(2):335--388, 1987.

\bibitem{Jones1994}
V.~F.~R. Jones.
\newblock The {P}otts model and the symmetric group.
\newblock In {\em Subfactors ({K}yuzeso, 1993)}, pages 259--267. World Sci.
  Publ., River Edge, NJ, 1994.

\bibitem{Jones1994_a}
V.~F.~R. Jones.
\newblock A quotient of the affine {H}ecke algebra in the {B}rauer algebra.
\newblock {\em Enseign. Math. (2)}, 40(3-4):313--344, 1994.

\bibitem{Kauffman1987}
L.~H. Kauffman.
\newblock State models and the {J}ones polynomial.
\newblock {\em Topology}, 26(3):395--407, 1987.

\bibitem{Kauffman1990}
L.~H. Kauffman.
\newblock An invariant of regular isotopy.
\newblock {\em Trans. Amer. Math. Soc.}, 318(2):417--471, 1990.

\bibitem{Kauffman1997}
L.~H. Kauffman.
\newblock Knots and diagrams.
\newblock In {\em Lectures at {KNOTS} '96 ({T}okyo)}, volume~15 of {\em Ser.
  Knots Everything}, pages 123--194. World Sci. Publ., River Edge, NJ, 1997.

\bibitem{KL1994}
L.~H. Kauffman and S.~L. Lins.
\newblock {\em Temperley-{L}ieb recoupling theory and invariants of
  {$3$}-manifolds}, volume 134 of {\em Annals of Mathematics Studies}.
\newblock Princeton University Press, Princeton, NJ, 1994.

\bibitem{LF2006}
K.~W. Lau and D.~G. FitzGerald.
\newblock Ideal structure of the {K}auffman and related monoids.
\newblock {\em Comm. Algebra}, 34(7):2617--2629, 2006.

\bibitem{NS1978}
T.~E. Nordahl and H.~E. Scheiblich.
\newblock Regular {$\ast $}-semigroups.
\newblock {\em Semigroup Forum}, 16(3):369--377, 1978.

\bibitem{RSA2014}
D.~Ridout and Y.~Saint-Aubin.
\newblock Standard modules, induction and the structure of the
  {T}emperley-{L}ieb algebra.
\newblock {\em Adv. Theor. Math. Phys.}, 18(5):957--1041, 2014.

\bibitem{TL1971}
H.~N.~V. Temperley and E.~H. Lieb.
\newblock Relations between the ``percolation'' and ``colouring'' problem and
  other graph-theoretical problems associated with regular planar lattices:
  some exact results for the ``percolation'' problem.
\newblock {\em Proc. Roy. Soc. London Ser. A}, 322(1549):251--280, 1971.

\bibitem{Westbury1995}
B.~W. Westbury.
\newblock The representation theory of the {T}emperley-{L}ieb algebras.
\newblock {\em Math. Z.}, 219(4):539--565, 1995.

\bibitem{Wilcox2007}
S.~Wilcox.
\newblock Cellularity of diagram algebras as twisted semigroup algebras.
\newblock {\em J. Algebra}, 309(1):10--31, 2007.

\end{thebibliography}
\bibliographystyle{abbrv}

\end{document}